\documentclass[3p,times,fleqn]{elsarticle}

\usepackage{hyperref}

\journal{Comput. Meth. Appl. Mech. Eng.}

\usepackage{amssymb,amsmath,mathtools}
\usepackage{amsthm}
\theoremstyle{remark}
\newtheorem{remark}{Remark}[section]
\theoremstyle{plain}
\newtheorem{theorem}{Theorem}[section]

\newtheorem{lemma}[theorem]{Lemma}

\usepackage{float}
\usepackage[section]{algorithm}
\usepackage{algpseudocode}
\usepackage{subfig}
\usepackage{multirow}
\newcommand*{\vcenteredhbox}[1]{\begingroup
\setbox0=\hbox{#1}\parbox{\wd0}{\box0}\endgroup}
\usepackage{color}
\usepackage[export]{adjustbox}
\usepackage{xspace}

\newcommand{\tensorTwo}[1]{\boldsymbol#1}
\newcommand{\tS}{\widetilde{S}}

\newcommand{\precname}{RACP\xspace}
\newcommand{\precnameVar}[1]{RACP(#1)\xspace}
\renewcommand{\vec}[1]{\mathbf{#1}}
\newcommand{\convec}[1]{\boldsymbol{#1}}
\newcommand{\Vecfunspace}[1]{\boldsymbol{\mathcal{#1}}}
\providecommand{\mat}[1]{#1}

\graphicspath{{figs/}}

\begin{document}

\begin{frontmatter}

\title{A Reverse Augmented Constraint preconditioner for Lagrange multiplier methods in contact mechanics}

\author[addr1]{Andrea Franceschini\corref{mycorrespondingauthor}}
\cortext[mycorrespondingauthor]{Corresponding author}
\ead{andrea.franceschini@unipd.it}

\author[addr2]{Matteo Frigo}
\author[addr2]{Carlo Janna}
\author[addr1]{Massimiliano Ferronato}

\address[addr1]{Department of Civil, Environmental and Architectural 
Engineering, University of Padova, Padova, Italy}
\address[addr2]{M3E S.r.l., via Giambellino 7, 35129 Padova, Italy}

\begin{abstract}
Frictional contact is one of the most challenging problems in computational mechanics. Typically, it is a tough non-linear problem often requiring several Newton iterations to converge and causing troubles also in the solution to the related linear systems. When contact is modeled with the aid of Lagrange multipliers, the impenetrability condition is enforced exactly, but the associated Jacobian matrix is indefinite and needs a special treatment for a fast numerical solution. In this work, a constraint preconditioner is proposed where the primal Schur complement is computed after augmenting the zero block. The name {\em Reverse} is used in contrast to the traditional approach where only the structural block undergoes an augmentation. Besides being able to address problems characterized by singular structural blocks, often arising in contact mechanics, it is shown that the proposed approach is significantly cheaper than traditional constraint preconditioning for this class of problems and it is suitable for an efficient HPC implementation through the Chronos parallel package.
Our conclusions are supported by several numerical experiments on mid- and large-size problems from various applications. The source files implementing the proposed algorithm are freely available on GitHub.

\end{abstract}


\end{frontmatter}

\section{Introduction}
\label{sec:intro}

Frictional contact mechanics is a very significant engineering problem that arises whenever an interaction exists between different bodies or parts of the same body. Motion of the entities in contact is usually activated at the occurrence of particular stress conditions and evolves under specific constraints, such as the impenetrability of solid bodies and the static-dynamic friction law. In this work, we refer to the equilibrium of a deformable solid occupying the domain $\Omega\subset\mathbb{R}^3$
with the assumptions of quasi-static conditions and infinitesimal strains. The standard linear momentum balance differential equation 
%
\begin{equation}
-\nabla \cdot \tensorTwo{\sigma} (\convec{u}) = \convec{b}
  \label{eq:momentumBalanceS}
\end{equation}
is subject to the following conditions over the contact surface $\Gamma_f$ with normal direction $\convec{n}_f$ \cite{KikOde88,Lau03,Wri06}:
\begin{subequations}
  \begin{align}
  & t_N = \convec{t} \cdot \convec{n}_f \le 0,
  &
  & g_N = \llbracket \convec{u} \rrbracket \cdot \convec{n}_f \ge 0,
  &
  & t_N g_N = 0,
  &
  & \mbox{(impenetrability)},
  \label{eq:normal_contact_KKT} \\
  & \left\| \convec{t}_T \right\|_2 \le \tau_{\max}(t_N),
  &
  & \dot{\convec{g}}_T \cdot \convec{t}_T = \tau_{\max}(t_N) || \dot{\convec{g}}_T ||_2,
  &&&
  & \mbox{(friction)}.
  \label{eq:frictional_contact_KKT}
\end{align}
\label{eq:KKT}\null
\end{subequations}
In the inequality-constrained problem \eqref{eq:momentumBalanceS}-\eqref{eq:KKT}, the displacement $\convec{u}$ is the primary unknown, with: $\convec{b}$ the external body forces; $\tensorTwo{\sigma}
(\convec{u})$ the stress tensor; $\convec{t} = t_N \convec{n}_f + \convec{t}_T$ the traction
over $\Gamma_f$, decomposed into its normal and tangential components, $t_N$ and $\convec{t}_T$; $\llbracket \convec{u} \rrbracket = 
g_N \convec{n}_f + \convec{g}_T$ the jump of $\convec{u}$ across $\Gamma_f$, decomposed into its normal and tangential components, $g_N$ and $\convec{g}_T$;
and $\tau_{\max}(t_N)$ a bounding value for the measure of $\convec{t}_T$. Relationships \eqref{eq:KKT} are the Karush-Kuhn-Tucker (KKT) complementary conditions for normal and frictional contact \cite{SimHug98}. In essence, they state that: (i) the normal traction must be compressive if the contact exists, with no penetration allowed between the two sides of the discontinuity surface $\Gamma_f$ (equation \eqref{eq:normal_contact_KKT}), and (ii) an upper bound for the magnitude of the tangential component of traction is set, at which slip is allowed and is collinear with friction (equation \eqref{eq:frictional_contact_KKT}). The mathematical problem is closed by adding the constitutive relationships for the stress $\tensorTwo{\sigma}(\convec{u})$ and friction $\tau_{\max}(t_N)$, and prescribing appropriate Dirichlet and Neumann boundary conditions.

A well-posed variational formulation of problem \eqref{eq:momentumBalanceS}-\eqref{eq:KKT} can be obtained by either penalty regularization or Lagrange multipliers. The former approach introduces into the variational principle associated to the equilibrium equation a high cost for the constraint violation, so as to satisfy approximately the contact conditions. Despite its inexact nature, the penalty method is still quite popular, especially in the world of finite elements for solid mechanics \cite{Zie00,Bat06}, with significant applications in different engineering fields, e.g. \cite{KhoLew99,OnaRoj04,FerGamJanTea08,BenEssFak16,BurErn17}. By distinction, the latter approach prescribes the minimization of the constrained variational principle in a mathematically exact way by using Lagrange multipliers \cite{Ber84,Wri06}. Convergence and numerical stability of the non linear problem is generally improved \cite{HagHueWoh08,fraferjantea16} at the cost of adding new variables as primary unknowns and increasing the overall problem size.

Here, we focus on the Lagrange multiplier approach. In the context of contact mechanics, Lagrange multipliers have the physical meaning of traction vector $\convec{t}$ living on the discontinuity surface $\Gamma_f$. Denoting by $\Vecfunspace{U}$ and $\Vecfunspace{U}_0$ the subspace of $[H^1(\Omega)]^3$ acting as trial and test spaces for the displacement, respectively, the function space $\Vecfunspace{M}(\convec{t})$ for the Lagrange multipliers can be defined as a subspace of the dual space $\Vecfunspace{M}$ of the trace $\Vecfunspace{U}^{\Gamma_f}_0$ of $\Vecfunspace{U}_0$ restricted to $\Gamma_f$ such that \cite{KikOde88}:
\begin{equation}
    \Vecfunspace{M} (\convec{t}) = \left\{ \convec{\mu}\in\Vecfunspace{M}, \mu_N \leq 0, \left( \convec{\mu}, \convec{w} \right)_{\Gamma_f} \leq \left(\tau_{\max}(t_N),\|\convec{w}_T\|_2\right)_{\Gamma_f}, \convec{w}\in\Vecfunspace{U}^{\Gamma_f}_0 \mbox{ with } w_N \geq 0 \right\},
    \label{eq:defMt}
\end{equation}
where the pedix $N$ and $T$ denotes the normal and tangential component of any vector function defined on $\Gamma_f$, respectively. The weak variational form of \eqref{eq:momentumBalanceS}-\eqref{eq:KKT} consists of finding $\{\convec{u},\convec{t}\}\in\Vecfunspace{U}\times\Vecfunspace{M}(\convec{t})$ such that:
\begin{subequations}
  \begin{align}
      \left(\nabla^s\convec{\eta},\tensorTwo{\sigma}\right)_{\Omega} + \left(\llbracket\convec{\eta}\rrbracket,\convec{t}\right)_{\Gamma_f} = \left(\convec{\eta},\convec{b}\right)_{\overline{\Omega}}, & & &\forall \convec{\eta}\in\Vecfunspace{U}_0, \label{eq:virtualwork} \\
      \left(t_N-\mu_N,g_N\right)_{\Gamma_f} + \left(\convec{t}_T-\convec{\mu}_T,\dot{\convec{g}}_T\right)_{\Gamma_f} \geq 0, & & &\forall \convec{\mu}\in\Vecfunspace{M}(\convec{t}), \label{eq:congruence}
  \end{align}
  \label{eq:weak_form}
\end{subequations}
where \eqref{eq:virtualwork} expresses the virtual work principle and \eqref{eq:congruence} the compatibility conditions for the contact surface. The variational inequality \eqref{eq:congruence} can be transformed into an equality by detecting the current contact operating mode of every point lying on $\Gamma_f$, for instance with the aid of an active-set algorithm. Discretization of the continuous problem \eqref{eq:weak_form} is finally carried out by replacing the mixed function space $\Vecfunspace{U}\times\Vecfunspace{M}(\convec{t})$ with the discrete subspace $\Vecfunspace{U}^h\times\Vecfunspace{M}^h(\convec{t}^h)$ associated to a conforming partition of the geometrical domain. Setting $n_u=\dim{\Vecfunspace{U}^h}$ and $n_t=\dim{\Vecfunspace{M}^h(\convec{t}^h)}$, the Jacobian $\mathcal{J}$ of the discrete version of \eqref{eq:weak_form} turns out to be a generalized saddle-point matrix of size $n_u+n_t$:
\begin{equation}
    \mathcal{J} = \left[ \begin{array}{cc} A & B_1 \\ B_2 & 0 \end{array} \right],
    \label{eq:nonsym_Jac}
\end{equation}
where $A\in\mathbb{R}^{n_u\times n_u}$ is the tangent stiffness matrix of the continuous body, and $B_1\in\mathbb{R}^{n_u\times n_t}$, $B_2\in\mathbb{R}^{n_t\times n_u}$ are the rectangular blocks coupling displacements with tractions.

The efficient solution of the sequence of linear problems with the matrix $\mathcal{J}$ is the object of the present study. Systems of equations with a similar algebraic structure arise very frequently in the applications and a significant amount of recent literature has been devoted to this issue \cite{benzi2005}. Important examples of saddle-point problems can be found, for instance, in constrained optimization \cite{li2001existence,Bergamaschi2007137,schoberl2007symmetric,Pearson2020}, mixed finite element approximations \cite{vassilevski1996prec,barrientos2002mixed,loghin2004analysis}, Stokes' and Biot's models \cite{olshanskii2006uniform,bergamaschi2008mixed,axelsson2012stable,cao2015relaxed,Pearson2018331,hong2020parameter}, just to cite a few. Of course, the key for guaranteeing a robust, efficient and scalable solution of indefinite saddle-point systems is the definition and design of an appropriate preconditioning strategy. It is well recognized that the most effective paradigm for building a good preconditioner for the matrix \eqref{eq:nonsym_Jac} relies on properly approximating the leading block $A$ and the related Schur complement $S=-B_2A^{-1}B_1$ with some $\tilde{A}$ and $\tilde{S}$, respectively \cite{Murphy20001969,Keller20001300,benzi2005}. Such an operation, however, can be strongly problem-dependent. For instance,
a straightforward idea relies on setting $\tilde{A}$ equal to the diagonal of $A$, thus allowing for an explicit computation of the approximate Schur complement as $\tilde{S}=-B_2\tilde{A}^{-1}B_1$. This idea proves effective if $A$ is diagonally dominant, e.g.~\cite{Bergamaschi2004149}, but often fails in many other circumstances. Alternatively, explicit approximate inverses can be used for either $\tilde{A}^{-1}$, such as in \cite{Bergamaschi20072647,Janna2012661}, or the computation of the decoupling factors $F=A^{-1}B_1$ and $G=B_2A^{-T}$, e.g.~\cite{ferronato2019general,Nardean2021}, so that the approximate Schur complement can be recast as $\tilde{S}=-GAF$.
Other possibilities involve problem-dependent strategies for the specific Schur complement at hand, including also physics-based strategies, e.g.~\cite{Silvester2001261,elman2002preconditioners,cao2003fast,Elman20061651,choi2015practical,Castelletto2016894}.
All these techniques, however, can become critical when the leading block is singular, with theoretical results available only in case of maximal rank deficiency \cite{Greif20191}.

In this work, we focus on a reverse approach, which introduces a non-zero contribution at the (2,2) block, inspired by Augmented Lagrangian techniques \cite{fortin1983augmented,bacuta2006unified,Benzi20062095,lee2007robust,brezzi2012mixed}. We discuss the properties of the proposed preconditioner and verify its performance in comparison with other methods currently available for the same problem. Finally, an implementation prone for large-size massively parallel simulations is described and tested, showing the algorithm efficiency and scalability.
The proposed algorithm has been implemented on the top of the Chronos package \cite{CHRONOS-webpage}, and the \precname sources are freely available on GitHub at \url{https://github.com/matteofrigo5/HPC\_ReverseAugmentedConstrained/}.

\section{Preconditioning framework}
\label{sec:prec}

The classical way to address the preconditioning of matrix $\mathcal{J}$ of equation \eqref{eq:nonsym_Jac} relies on exploiting the inverse of its block LDU factorization:
\begin{equation}
    \mathcal{J}^{-1} = \left[ \begin{array}{cc} I_u & -A^{-1}B_1 \\ 0 & I_t \end{array} \right]
    \left[ \begin{array}{cc} A^{-1} & 0 \\ 0 & S^{-1} \end{array} \right]
    \left[ \begin{array}{cc} I_u & 0 \\ -B_2A^{-1} & I_t \end{array} \right],
    \label{eq:LDU_J}
\end{equation}
where $I_u$ and $I_t$ are the identity matrices of order $n_u$ and $n_t$, resepctively. In the factorization \eqref{eq:LDU_J},
both $A$ and the Schur complement $S=-B_2A^{-1}B_1$ are replaced by sparse approximations, $\tilde{A}$ and $\tilde{S}$, whose inverse should be easy and effective to apply.
While several off-the-shelf algebraic tools for $\tilde{A}^{-1}$ are available, such as incomplete factorizations, approximate inverses and multigrid methods, the same is not true for $\tilde{S}^{-1}$. First, it is troublesome to obtain a good sparse approximation of $S$, since the Schur complement is almost completely dense because of $A^{-1}$. Second, the inexact application of $\tilde{S}^{-1}$ introduces an additional level of approximation, which further contributes to the difficulty in designing an effective preconditioner.
Franceschini et al.~\cite{franceschini2019block} investigated the use of different explicit and implicit choices for the Schur complement matrix in contact mechanics problems, including a physics-based block diagonal approximation, a generalized least-square commutator and an explicit approximate inverse. All these methods, however, require the assumption for $A$ to be non singular, but this condition is not satisfied in several applications.

%

In this work, we focus on the solution of a linear system with
the symmetric saddle-point matrix:
\begin{equation}
  \mathcal{A} = \begin{bmatrix}
    A   & B \\
    B^T & 0
  \end{bmatrix},
  \label{eq:matA}
\end{equation}
with $A \in \mathbb{R}^{n_u \times n_u}$ a symmetric positive semidefinite (SPSD) matrix and $B \in
\mathbb{R}^{n_u \times n_t}$, $n_u>n_t$, a full rank constraint matrix. In order to accelerate the convergence of a Krylov supspace solver, we use as a preconditioner the inverse of the ``stabilized'' saddle-point matrix
\begin{equation}
    \hat{\mathcal{A}} = \begin{bmatrix}
    A   & B \\
    B^T & -C
  \end{bmatrix},
  \label{eq:matAhat}
\end{equation}
with $C\in\mathbb{R}^{n_t\times n_t}$ symmetric positive definite (SPD).

\begin{remark}
\label{rem:nonsym}
The saddle-point matrix $\mathcal{A}$ is the symmetrized version of $\mathcal{J}$ in equation \eqref{eq:nonsym_Jac}. Actually, the results presented here can be extended in practice to $\mathcal{J}$ as well, as it will be shown in the numerical experiments, even though some theoretical analyses are no longer valid.
\end{remark}

\begin{theorem}
\label{thm:ideal_eig}
Let $\mathcal{A}$ and $\hat{\mathcal{A}}$ be the saddle-point matrices defined in equation \eqref{eq:matA} and \eqref{eq:matAhat}, with $A$ non singular. If $C=B^TA^{-1}B$, then the eigenvalues $\lambda$ of the preconditioned matrix $\hat{\mathcal{A}}^{-1}\mathcal{A}$ are either 1, with multiplicity $n_u$, or 0.5, with multiplicity $n_t$.
\end{theorem}
\begin{proof}
The eigenvalues of $\hat{\mathcal{A}}^{-1}\mathcal{A}$ satisfy the generalized eigenproblem $\mathcal{A}\vec{v}=\lambda\hat{\mathcal{A}}\vec{v}$:
\begin{equation}
  \begin{bmatrix}
    A & B \\
    B^T & 0 \\
  \end{bmatrix} \begin{bmatrix}
    \vec{v}_u \\
    \vec{v}_t \\
  \end{bmatrix} = \lambda \begin{bmatrix}
    A & B \\
    B^T & -C \\
  \end{bmatrix} \begin{bmatrix}
    \vec{v}_u \\
    \vec{v}_t \\
  \end{bmatrix},
\end{equation}
that is
\begin{equation}
  \begin{cases}
    A\vec{v}_u+B\vec{v}_t = \lambda A\vec{v}_u + \lambda B\vec{v}_t \\
    B^T \vec{v}_u = \lambda B^T \vec{v}_u - \lambda C\vec{v}_t \\
  \end{cases}
  \label{eq:geneig}
\end{equation}
The first equation
\begin{equation}
  \left(1-\lambda\right)\left(A\vec{v}_u+B\vec{v}_t\right) = \vec{0}_u,
\end{equation}
yields either $\lambda = 1$ or $A\vec{v}_u = -B\vec{v}_t$. If $\lambda=1$, it follows from the second equation in \eqref{eq:geneig} that $C\vec{v}_t = \vec{0}_t$, which implies $\vec{v}_t = \vec{0}_t$ because of the regularity of $C$, while $\vec{v}_u$ can be any non-zero vector in $\mathbb{R}^{n_u}$. 
Therefore, $\lambda=1$ has multiplicity $n_u$.
If $\lambda \ne 1$, then $\vec{v}_u = -A^{-1}B\vec{v}_t$ and 
the second equation in \eqref{eq:geneig} becomes:
\begin{equation}
  -B^T A^{-1} B \vec{v}_t = \lambda \left(-C - B^T A^{-1} B\right) \vec{v}_t.
  \label{eq:Ceig}
\end{equation}
Recalling the assumption on $C$,
we have:
\begin{equation}
  \left(C\right)^{-1} C \vec{v}_t = 2\lambda \vec{v}_t,
\end{equation}
which implies $\lambda=0.5$ for any non-zero vector $\vec{v}_t$ in $\mathbb{R}^{n_t}$.
\end{proof}

The outcome of Theorem \ref{thm:ideal_eig} suggests the ideal choice for $C$ if the leading block $A$ is regular. If not, we can still use $\hat{\mathcal{A}}^{-1}$ as a preconditioner by setting $C\neq B^TA^{-1}B$ and using its block UDL factorization:
\begin{equation}
  \hat{\mathcal{A}}=\mathcal{U}\mathcal{D}\mathcal{L} = \begin{bmatrix}
    I_u & -B C^{-1} \\
    0 & I_t \\
  \end{bmatrix} \begin{bmatrix}
    S_u & 0 \\
    0 & -C \\
  \end{bmatrix} \begin{bmatrix}
    I_u &  0 \\
    -C^{-1} B^T & I_t \\
  \end{bmatrix},
  \label{eq:bUDL_sym}
\end{equation}
with $S_u$ the so-called ``primal'' Schur complement \cite{benzi2005} defined as $S_u = A + B C^{-1} B^T$. Note that $S_u \in
\mathbb{R}^{n_u \times n_u}$ is SPD.
Since $n_u$ is usually large, the inverse of $S_u$ has to be applied inexactly through some appropriate operator $\tS_u^{-1}$.
The resulting preconditioner
$\mathcal{M}^{-1}$ finally reads:
\begin{equation}
  \mathcal{M}^{-1} = \mathcal{L}^{-1}\mathcal{\widetilde{D}}^{-1}\mathcal{U}^{-1} =
  \begin{bmatrix}
    I_u &  0 \\
    C^{-1} B^T & I_t \\
  \end{bmatrix} \begin{bmatrix}
    \tS_u^{-1} & 0 \\
    0 & -C^{-1} \\
  \end{bmatrix} \begin{bmatrix}
    I_u & B C^{-1} \\
    0 & I_t \\
  \end{bmatrix}.
  \label{eq:RACPsym}
\end{equation}
In essence, $\mathcal{M}^{-1}$ uses a reverse approach with respect to classical constraint preconditioners, by exploiting the primal Schur complement of an augmented matrix. For this reason, we denote it as Reverse Augmented Constraint Preconditioner (\precname).

\subsection{Eigenvalue analysis}
\label{sec:bounds}

The eigenvalues of the \precname preconditioned matrix $\mathcal{M}^{-1}\mathcal{A}$ satisfy the bounds that follow.

\begin{theorem}
\label{thm:sym_bounds}
Let $\mathcal{A}$ and $\mathcal{M}^{-1}$ be the matrices defined in equation \eqref{eq:matA} and \eqref{eq:RACPsym}, respectively, and:
\begin{subequations}
\begin{align}
  &\alpha_u = 
  \lambda_{\min}
  \left(\tS_u^{-1}
    \left(S_u + B C^{-1} B^T\right)\right), \qquad 
  \beta_u = 
  \lambda_{\max}
  \left(\tS_u^{-1}
    \left(S_u + B C^{-1} B^T\right)\right), 
    \label{eq:au_bu} \\
  &\alpha_t = 
  \sigma_{\min}
    \left( \tS_u^{-1/2}BC^{-1/2} \right), \qquad 
  \beta_t = 
  \sigma_{\max}
    \left( \tS_u^{-1/2}BC^{-1/2} \right). \label{eq:at_bt}
\end{align}
\label{eq:a_b}
\end{subequations}
where $\lambda(\cdot)$ and $\sigma(\cdot)$ denote eigenvalues and singular values,
respectively, of the matrix within brackets.
Then, the real eigenvalues of $\mathcal{M}^{-1}\mathcal{A}$ are such that:
\begin{equation}
  \min\left\{\alpha_u,\frac{2 \alpha_t^2}{\beta_u + \sqrt{\beta_u^2 - 4 \alpha_t^2}}
    \right\} \le \lambda \le \beta_u,
    \label{eq:rebound_sym}
\end{equation}
and the real and imaginary part, $\lambda_{\Re}$ and $\lambda_{\Im}$, of the complex eigenvalues are such that:
\begin{equation}
  \frac{\alpha_u}{2} \le \lambda_{\Re} \le \frac{\beta_u}{2}, \qquad 
    |\lambda_{\Im}| \le \sqrt{\beta_t^2 - \frac{\alpha_u^2}{4}},
    \label{eq:imbound_sym}
\end{equation}
with no complex eigenvalues if $2\beta_t < \alpha_u$.
\end{theorem}
\begin{proof}
The eigenvalues of $\mathcal{M}^{-1}\mathcal{A}$
are the solution of the generalized eigenproblem:
\begin{equation}
  \mathcal{A} \vec{v} = \lambda \mathcal{M} \vec{v}.
  \label{eq:geneig_sym}
\end{equation}
Recalling the factorization \eqref{eq:RACPsym} and setting $\vec{w}=\mathcal{L}\vec{v}$, we have:
\begin{equation}
  \begin{bmatrix}
    S_u + B C^{-1} B^T & B \\
    B^T & 0 \\
  \end{bmatrix} \vec{w} = \lambda \begin{bmatrix}
    \tS_u & 0 \\
    0 & -C \\
  \end{bmatrix} \vec{w}, 
  \label{eq:geneig_sym2}
\end{equation}
which is equivalent to:
\begin{equation}
  \begin{bmatrix}
    \tS_u^{-\frac{1}{2}} \left(S_u + B C^{-1} B^T\right) \tS_u^{-\frac{1}{2}} &
      \tS_u^{-\frac{1}{2}} B C^{-\frac{1}{2}} \\
    -C^{-\frac{1}{2}} B^T \tS_u^{-\frac{1}{2}} & 0
  \end{bmatrix} \vec{z} = 
  \lambda \vec{z}, \qquad
  \vec{z} = \begin{bmatrix}
    \tS_u^{\frac{1}{2}} & 0 \\
    0 & C^{\frac{1}{2}}
  \end{bmatrix} \vec{w}.
  \label{eq:geneig_sym3}
\end{equation}
Hence, $\mathcal{M}^{-1}\mathcal{A}$ is similar to a non-symmetric saddle point matrix, for which
the eigenvalue bounds
\eqref{eq:rebound_sym} and \eqref{eq:imbound_sym} hold
\cite{bergamaschi2012eigenvalue}.
\end{proof}


\begin{lemma}
\label{lem:sym_bounds}
In the limiting case $\tS_u^{-1}=S_u^{-1}$, the bounds \eqref{eq:rebound_sym}-\eqref{eq:imbound_sym} read:
\begin{subequations}
\begin{align}
  &\min\left\{1,\frac{2 \alpha_t^2}{1+\beta_t^2 + \sqrt{(1+\beta_t^2)^2 - 4 \alpha_t^2}}
    \right\} \le \lambda \le 1+\beta_t^2,
    \label{eq:rebound_symex} \\
  &\frac{1}{2} \le \lambda_{\Re} \le \frac{1+\beta_t^2}{2}, \qquad 
    |\lambda_{\Im}| \le \sqrt{\beta_t^2 - \frac{1}{4}},
    \label{eq:imbound_symex}
    \end{align}
    \label{eq:bounds_symex}
\end{subequations}
with no complex eigenvalues if $\beta_t<1/2$.
\end{lemma}
\begin{proof}
The leading block $
\tS_u^{-\frac{1}{2}} \left(S_u + B C^{-1} B^T\right) \tS_u^{-\frac{1}{2}}$ of the matrix in equation \eqref{eq:geneig_sym3} with $\tS_u^{-1}=S_u^{-1}$ 
becomes
    $
    I_u + F F^T$,
with $F=\tS_u^{-\frac{1}{2}} B C^{-\frac{1}{2}}$. Since $F$ has rank $n_t<n_u$,
we have $\alpha_u = 1$ and $\beta_u =
\beta_t^2 + 1$.
Then, the result \eqref{eq:bounds_symex} follows immediately from \eqref{eq:rebound_sym}-\eqref{eq:imbound_sym}.
\end{proof}

The bounds provided by Lemma \ref{lem:sym_bounds} represent the limiting condition that can be potentially achieved with the \precname algorithm \eqref{eq:RACPsym} for any SPD matrix $C$. Despite both $\mathcal{M}^{-1}$ and $\mathcal{A}$ are symmetric, the preconditioned matrix $\mathcal{M}^{-1}\mathcal{A}$ is not and complex eigenvalues might arise. Symmetry can be restored if we define the \precname operator starting from the matrix:
\begin{equation}
  \overline{\mathcal{A}} = \begin{bmatrix}
    A & B \\
    -B^T & C \\
  \end{bmatrix} 
= \begin{bmatrix}
    I_u & B C^{-1} \\
    0 & I_t \\
  \end{bmatrix} \begin{bmatrix}
    S_u & 0 \\
    0 & C \\
  \end{bmatrix} \begin{bmatrix}
    I_u & 0 \\
    -C^{-1} B^T & I_t
  \end{bmatrix}=\mathcal{U}_a \mathcal{D}_a \mathcal{L}_a.
  \label{eq:matAbar}
\end{equation}
Approximating the
inverse of $S_u$ with $\tS_u^{-1}$, we can build an alternative \precname operator $\mathcal{M}_a^{-1}$ as:
\begin{equation}
  \mathcal{M}_a^{-1} = \mathcal{L}_a^{-1}\mathcal{\widetilde{D}}_a^{-1}\mathcal{U}_a^{-1} =
  \begin{bmatrix}
    I_u &  0 \\
    C^{-1} B^T & I_t \\
  \end{bmatrix} \begin{bmatrix}
    \tS_u^{-1} & 0 \\
    0 & C^{-1} \\
  \end{bmatrix} \begin{bmatrix}
    I_u & -B C^{-1} \\
    0 & I_t \\
  \end{bmatrix}.
  \label{eq:RACPnonsym}
\end{equation}
For the eigenvalues of the preconditioned matrix $\mathcal{M}_a^{-1}\mathcal{A}$ the following result holds.

\begin{theorem}
\label{thm:nonsym_bounds}
Let $\mathcal{A}$ and $\mathcal{M}_a^{-1}$ be the matrices defined in equation \eqref{eq:matA} and \eqref{eq:RACPnonsym}, respectively, and:
\begin{equation}
  \alpha_a = 
  \lambda_{\min}
  \left(\tS_u^{-1} A\right),
    \quad 
  \beta_a = 
  \lambda_{\max}
  \left(\tS_u^{-1} A\right),
    \quad
  \alpha_t = 
  \sigma_{\min}
    \left( \tS_u^{-1/2}BC^{-1/2} \right), \quad 
  \beta_t = 
  \sigma_{\max}
    \left( \tS_u^{-1/2}BC^{-1/2} \right), 
\label{eq:a_b_nonsym}
\end{equation}
where $\lambda(\cdot)$ and $\sigma(\cdot)$ denote eigenvalues and singular values,
respectively, of the matrix within brackets.
Then, the eigenvalues of $\mathcal{M}_a^{-1}\mathcal{A}$ are all real and such that:
\begin{equation}
  \lambda \in \left[\frac{\alpha_a - \sqrt{\alpha_a^2 + 4 \beta_t^2}}{2},
    \frac{\beta_a - \sqrt{\beta_a^2 + 4 \alpha_t^2}}{2}\right] \cup \left[\alpha_a,
    \frac{\beta_a + \sqrt{\beta_a^2 + 4 \beta_t^2}}{2}\right].
    \label{eq:bound_nonsym}
\end{equation}
\end{theorem}
\begin{proof}
To compute the eigenvalues of $\mathcal{M}_a^{-1}\mathcal{A}$, we observe that it is similar
to the symmetric saddle-point matrix:
\begin{equation}
\mathcal{\widetilde{D}}_a^{-\frac{1}{2}} \mathcal{U}_a^{-1} \mathcal{A} \mathcal{L}_a^{-1} \mathcal{\widetilde{D}}_a^{-\frac{1}{2}} =
  \begin{bmatrix}
    \tS_u^{-\frac{1}{2}} A \tS_u^{-\frac{1}{2}} & \tS_u^{-\frac{1}{2}} B C^{-\frac{1}{2}} \\
    C^{-\frac{1}{2}} B^T \tS_u^{-\frac{1}{2}} & 0 \\
  \end{bmatrix}, 
\label{eq:sim_nonsym}
\end{equation}
for which the bound \eqref{eq:bound_nonsym} holds
\cite{axelsson2006eigenvalue,ruiz2018refining}. 
%
\end{proof}

\begin{remark}
\label{rem:comp_bound}
The result of the eigenvalue analysis provided by Theorem \ref{thm:sym_bounds} and \ref{thm:nonsym_bounds} shows that the \precname formulation \eqref{eq:RACPsym} produces a non-symmetric preconditioned matrix with possibly complex eigenvalues. The real part $\lambda_{\Re}$ is always positive, while the imaginary part $\lambda_{\Im}$ satisfies a narrower bound than $\lambda_{\Re}$. By distinction, the alternative \precname formulation \eqref{eq:RACPnonsym} produces a symmetric but indefinite preconditioned matrix, with eigenvalues possibly lying in nearly opposite intervals with respect to the origin.
\end{remark}

\begin{theorem}
\label{thm:ideal_eig_a}
Let $\mathcal{A}$ and $\overline{\mathcal{A}}$ be the saddle-point matrices defined in
equation \eqref{eq:matA} and \eqref{eq:matAbar}, with $A$ non singular. Then the
eigenvalues $\lambda_a$ of the preconditioned matrix $\overline{\mathcal{A}}^{-1}
\mathcal{A}$ are either 1, with multiplicity $n_u$, or the opposite of the non-unitary
eigenvalues of $\hat{\mathcal{A}}^{-1} \mathcal{A}$. In particular, if $C=B^TA^{-1}B$,
then the $n_t$ non-unitary eigenvalues are $-0.5$.
\end{theorem}
\begin{proof}
The eigenvalues of $\overline{\mathcal{A}}^{-1}\mathcal{A}$ satisfy the generalized
eigenproblem $\mathcal{A}\vec{v}=\lambda_a\overline{\mathcal{A}}\vec{v}$:
\begin{equation}
  \begin{bmatrix}
    A & B \\
    B^T & 0 \\
  \end{bmatrix} \begin{bmatrix}
    \vec{v}_u \\
    \vec{v}_t \\
  \end{bmatrix} = \lambda_a \begin{bmatrix}
    A & B \\
    -B^T & C \\
  \end{bmatrix} \begin{bmatrix}
    \vec{v}_u \\
    \vec{v}_t \\
  \end{bmatrix},
\end{equation}
that is
\begin{equation}
  \begin{cases}
    A\vec{v}_u+B\vec{v}_t = \lambda_a A\vec{v}_u + \lambda_a B\vec{v}_t \\
    B^T \vec{v}_u = -\lambda_a B^T \vec{v}_u + \lambda_a C\vec{v}_t \\
  \end{cases}
  \label{eq:geneig_a}
\end{equation}
As before, the first equation is:
\begin{equation}
  \left(1-\lambda_a\right)\left(A\vec{v}_u+B\vec{v}_t\right) = \vec{0}_u,
\end{equation}
and it yields either $\lambda_a = 1$ or $A\vec{v}_u = -B\vec{v}_t$. If $\lambda_a=1$, it
follows from the second equation in \eqref{eq:geneig_a} that $C\vec{v}_t = 2 B^T
\vec{v}_u$, which implies a relation between $\vec{v}_u$ and $\vec{v}_t$ that, given the
regularity of $C$, can always be satisfied. Therefore, $\lambda_a=1$ has multiplicity
$n_u$. If $\lambda_a \ne 1$, then $\vec{v}_u = -A^{-1}B\vec{v}_t$ and the second equation
in \eqref{eq:geneig_a} becomes:
\begin{align}
  -B^T A^{-1} B \vec{v}_t &= \lambda_a B^T A^{-1} B \vec{v}_t + \lambda_a C \vec{v}_t,
    \nonumber\\
  -B^T A^{-1} B \vec{v}_t &= \frac{\lambda_a}{1+\lambda_a} C \vec{v}_t.
  \label{eq:Ceig_a}
\end{align}
Thus, the eigenvalues of $\overline{\mathcal{A}}^{-1}\mathcal{A}$ satisfy the relation:
\begin{equation}
  \frac{\lambda_a}{1+\lambda_a} = - \sigma\left[C^{-1} \left(B^T A^{-1} B\right)
    \right] = -\tau \quad \Rightarrow \quad \lambda_a = -\frac{\tau}{1+\tau}.
  \label{eq:lam_tau}
\end{equation}
From Theorem \ref{thm:ideal_eig}, the non-unitary eigenvalues of $\hat{\mathcal{A}}^{-1}
\mathcal{A}$ are $\sigma\left[\left(C + B^T A^{-1} B\right)^{-1}\left(B^T A^{-1} B\right)
\right]$. After some algebra, we have:
\begin{align}
  \left(C + B^T A^{-1} B\right)^{-1}\left(B^T A^{-1} B\right) &=
  \left[\left(B^T A^{-1} B\right)^{-1}\left(C + B^T A^{-1} B\right)\right]^{-1} =
    \nonumber \\
    &= \left[I + \left(B^T A^{-1} B\right)^{-1}C\right]^{-1} =
    \left\{I + \left[C^{-1}\left(B^T A^{-1} B\right)\right]^{-1}\right\}^{-1},
\end{align}
and, recalling that $\tau = \sigma\left[C^{-1} \left(B^T A^{-1} B\right)\right]$, we can
write:
\begin{equation}
  \lambda = \sigma\left[\left(C + B^T A^{-1} B\right)^{-1}\left(B^T A^{-1} B\right)\right]
    = \frac{\tau}{1 + \tau},
\end{equation}
thus, from Eq. \eqref{eq:lam_tau}, we have that $\lambda = -\lambda_a$.
\end{proof}

\subsection{Selection of matrix C}
\label{sec:matC}

Of course, one of the key factors for the effectiveness of the \precname approach is the choice of matrix $C$.
According to Theorem \ref{thm:ideal_eig}, $C$ should be close to the Schur complement $S=B^TA^{-1}B$ in order to cluster the eigenvalues of the preconditioned matrix around two non-zero values only. Since $C^{-1}$ is also needed to build explicitly the primal Schur complement $S_u$, a convenient option consists of defining $C$ as a diagonal matrix spectrally equivalent to $B^TA^{-1}B$.

The same problem was already addressed in \cite{franceschini2019block} by exploiting physics-based assumptions related to the locality of deformation. Since a variation of stress around a fracture is expected to produce a relevant displacement only in the close vicinity of the fracture itself, the full solution to the mechanical problem $A^{-1}\mathbf{b}_i$, with $\mathbf{b}_i$ the $i$-th column of $B$, $i=1,\ldots,n_t$, is replaced by a sparse vector with non-zero entries located only at the positions of the non-zeros in $\mathbf{b}_i$. The diagonal entries of $S$, which can be retained to build $C$, can be therefore approximated as:
\begin{equation}
  C_{i,i} = S_{i,i} \simeq r(\mathbf{b}_i^T) A\vert_{b_i}^{-1} r(\mathbf{b}_i) \qquad \forall i \in 1,\ldots,n_t,
  \label{eq:Cdef_fs}
\end{equation}
where $A\vert_{b_i}$ is a square block gathered from matrix $A$ at the row and column locations intercepting the non-zero entries in $\mathbf{b}_i$ and $r(\cdot)$ is a restriction operator on the vector within brackets retaining the non-zero entries only. This approximation was also used in \cite{Aagaard20133059,Jha20143776} for computing the increment in Newton's algorithm, or introduced in \cite{White201655,Castelletto2016894} for preconditioning purposes in coupled poromechanics.

Equation \eqref{eq:Cdef_fs} cannot be safely used if $A$ is rank-deficient, because $A\vert_{b_i}$ could be singular. Instead, we can refer to the classical methods used with augmented Lagrangian approaches \cite{Pow69,Golub20032076,benzi2005}, where the objective is to force the terms of the contribution $BC^{-1}B^T$ in $S_u$ to be of the same magnitude as those of $A$. The simplest way consists of setting $C=\gamma I_t$, with $\gamma=\|B\|^2/\|A\|$ for some convenient matrix norm \cite{Golub20032076,benzi2005}. However, using such a $C$ matrix, which is based on a global value of $\gamma$, may not represent a good choice as it may worsen the eigenspectrum and condition number of $S_u$ with respect to $A$. Instead, we use a local variant of this augmentation matrix that takes into account the physical observations leading to \eqref{eq:Cdef_fs}:
\begin{equation}
  C_{i,i} = \frac{\omega ||r(\mathbf{b}_i)||_2^2}{||A\vert_{b_i}||_2} \qquad \forall i \in 1,\ldots,n_t,
  \label{eq:Cdef}
\end{equation}
where
$\omega>0$ is a user-specified real relaxation parameter. The preconditioning algorithm obtained with the matrix $C$ defined as in \eqref{eq:Cdef} is denoted as \precnameVar{$\omega$}.

\begin{remark}
\label{rem:spar_pat}
Another important issue in the proposed approach is the sparsity pattern of the primal Schur complement $S_u$. Of course, it is desirable to change at the least possible the pattern of $A$, avoiding an excessive fill-in. Choosing a diagonal matrix as $C$, the pattern of $S_u$ differs from that of $A$ by $BB^T$. The matrix $B$ provides the connections between the degrees of freedom associated to the nodes lying on both sides of each fracture. Therefore, the sparsity pattern of $S_u$ corresponds to that of the elasticity matrix associated to a 3D body where the fractures are filled by elements with some fictitious size and stiffness. This guarantees that $S_u$ preserves a workable sparsity with no significant fill-in.
\end{remark}

\begin{figure}
  \centering
  \subfloat{\includegraphics[height=0.33\linewidth]{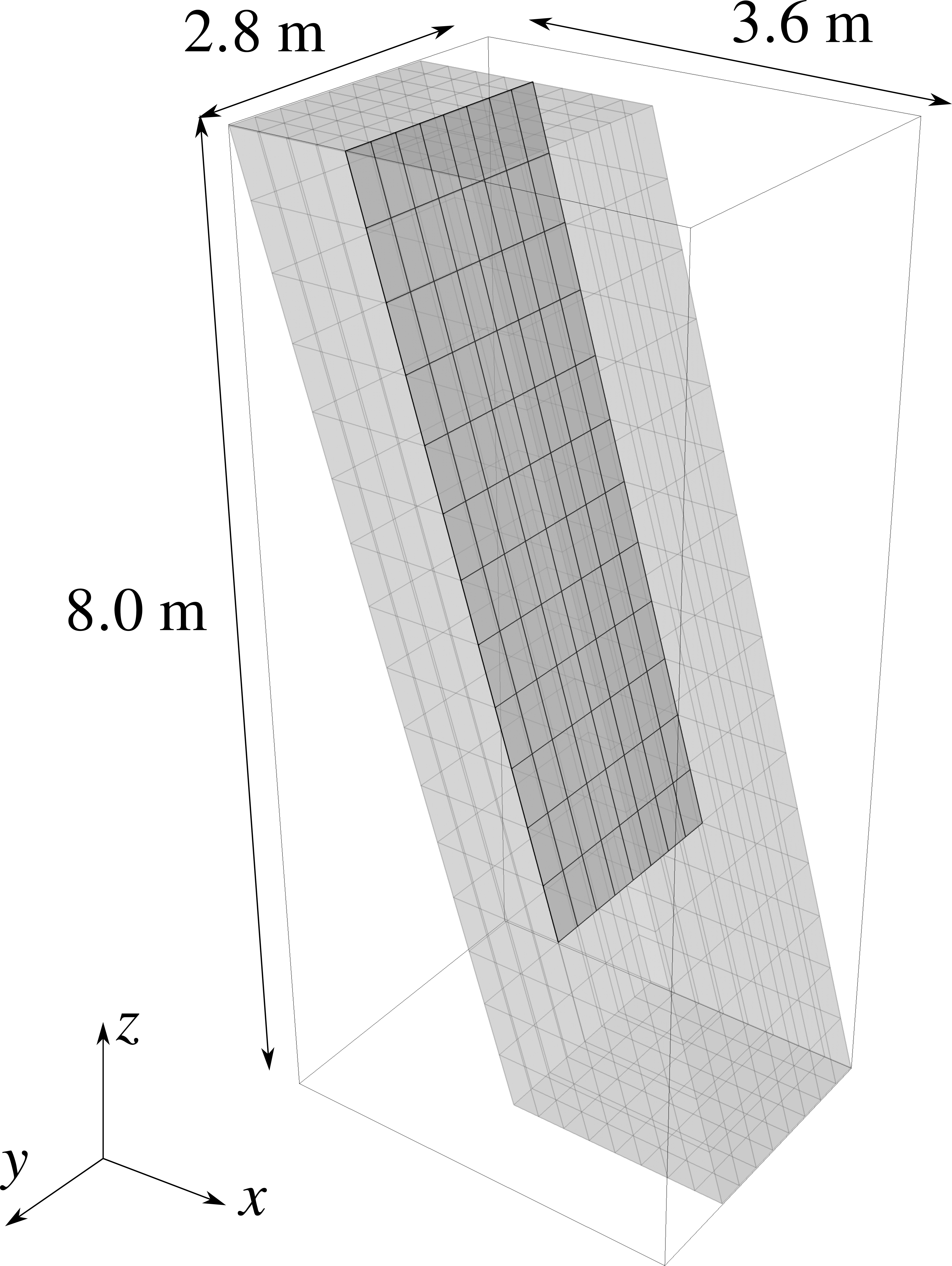}}
  \caption{Structured distorted grid with a single vertical fracture.} 
  \label{fig:cube}
\end{figure}

The effectiveness of the choices above is verified in the simple example shown in Figure \ref{fig:cube}, consisting of a structured distorted grid with a single fracture. The size of the problem blocks is $n_u=4455$ and $n_t=324$. First, we set $\widetilde{S}_u^{-1}=S_u^{-1}$ and check the eigenvalue bounds provided by Lemma \ref{lem:sym_bounds} and Theorem \ref{thm:nonsym_bounds} for $\mathcal{M}^{-1}\mathcal{A}$ and $\mathcal{M}_a^{-1}\mathcal{A}$, respectively. The bounds have been computed by setting $C$ as a diagonal matrix with entries equal to the values in either equation \eqref{eq:Cdef_fs} or \eqref{eq:Cdef} with $\omega=1$.
The results reported in Table \ref{tab:bounds} show that there is no significant difference between the two approaches for defining $C$.

\begin{table}
  \centering
  \begin{tabular}{c|c|c|c|c|c|c|c|c}
     & \multicolumn{4}{c|}{$\mathcal{M}^{-1}\mathcal{A}$ (Lemma \ref{lem:sym_bounds})} & \multicolumn{4}{c}{$\mathcal{M}_a^{-1}\mathcal{A}$ (Theorem \ref{thm:nonsym_bounds})} \\
     & \multicolumn{2}{c|}{$\lambda$} & $\lambda_{\Re}^*$ & $\lambda_{\Im}$ &
      \multicolumn{2}{c|}{$\lambda<0$} &
      \multicolumn{2}{c}{$\lambda>0$} \\
    matrix $C$ & lower & upper & upper & abs value & lower & upper & lower & upper \\
    \hline
    Eq. \eqref{eq:Cdef_fs} & 1.76e-1 & 2.00e+0 & 1.00e+0 & 8.66e-1 & -1.25e+0 &
     -2.58e-1 & 7.20e-4 & 1.85e+0 \\
    Eq. \eqref{eq:Cdef} & 1.78e-1 & 2.00e+0 & 1.00e+0 & 8.66e-1 & -9.99e-1 &
     -2.58e-1 & 6.88e-4 & 1.62e+0 \\
  \end{tabular}
  \caption{Theoretical eigenvalue bounds for $\mathcal{M}^{-1}\mathcal{A}$ and $\mathcal{M}^{-1}_a\mathcal{A}$ with $\widetilde{S}_u^{-1}=S_u^{-1}$ for the test case of Figure \ref{fig:cube}. 
    $^*$ according to equation \eqref{eq:imbound_symex}, the lower bound for the real part
    of any complex eigenvalue is 0.5, thus it has been omitted.}
  \label{tab:bounds}
\end{table}

\begin{figure}
  \centering
  \includegraphics[width=0.49\linewidth]{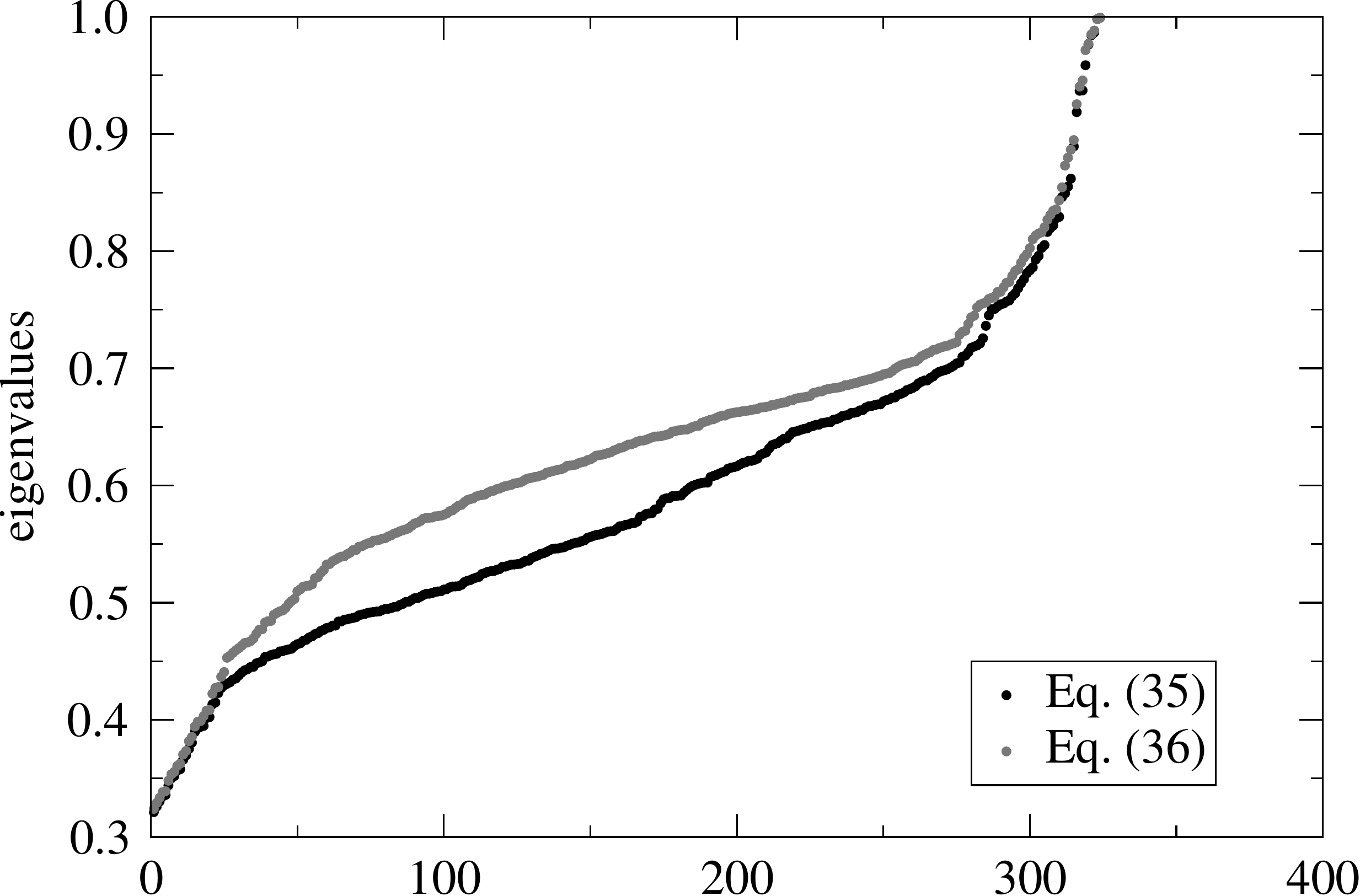}
  \caption{Non-unitary eigenvalues distribution of $\mathcal{M}^{-1}\mathcal{A}$ and $-\mathcal{M}_a^{-1}\mathcal{A}$ with $\tS_u^{-1}=S_u^{-1}$ for the test case of Figure \ref{fig:cube}.} 
  \label{fig:eigsPinvA}
\end{figure}

\begin{remark}
\label{rem:eig_symvsnosym}
It is easy to observe that setting $\widetilde{S}_u^{-1}=S_u^{-1}$ yields $n_u$ unitary eigenvalues for both $\mathcal{M}^{-1}\mathcal{A}$ and $\mathcal{M}_a^{-1}\mathcal{A}$. As proved by Theorem \ref{thm:ideal_eig_a}, the remaining $n_t$ eigenvalues are the same as those of equation \eqref{eq:Ceig} with either the positive or the negative sign at the left-hand side. Therefore, with $\widetilde{S}_u^{-1}=S_u^{-1}$ the non-unitary eigenvalues of the preconditioned matrix are all real positive for $\mathcal{M}^{-1}\mathcal{A}$ and negative for $\mathcal{M}_a^{-1}\mathcal{A}$. Their distribution in absolute value is provided in Figure \ref{fig:eigsPinvA} for the two choices of matrix $C$, showing that
also the full eigenspectrum does not appear to be significantly affected by either selection.
\end{remark}

\begin{figure}
  \centering
  \null\hfill
  \subfloat{\includegraphics[width=0.49\linewidth,valign=t]{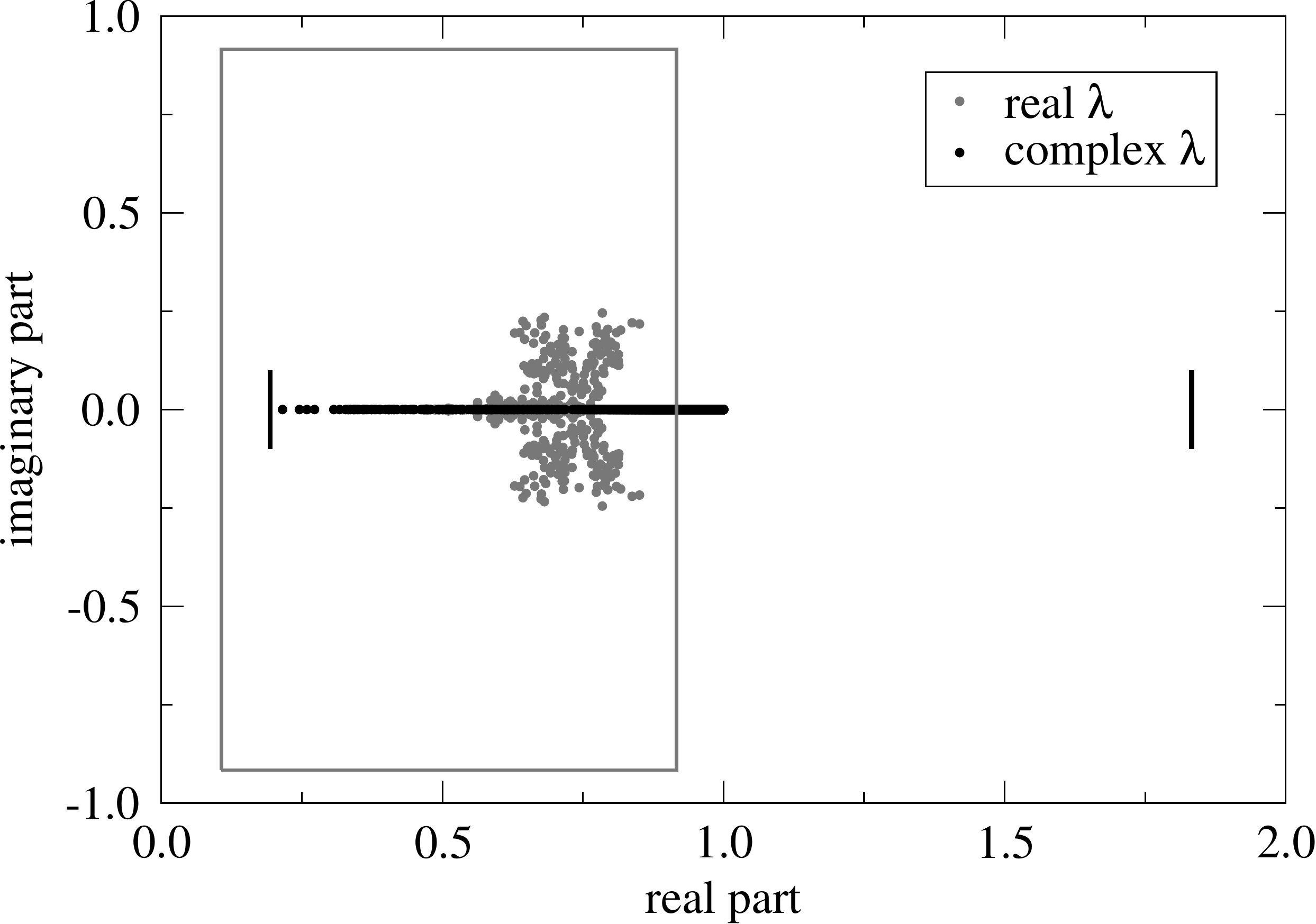}}\hfill
  \subfloat{\includegraphics[width=0.49\linewidth,valign=t]{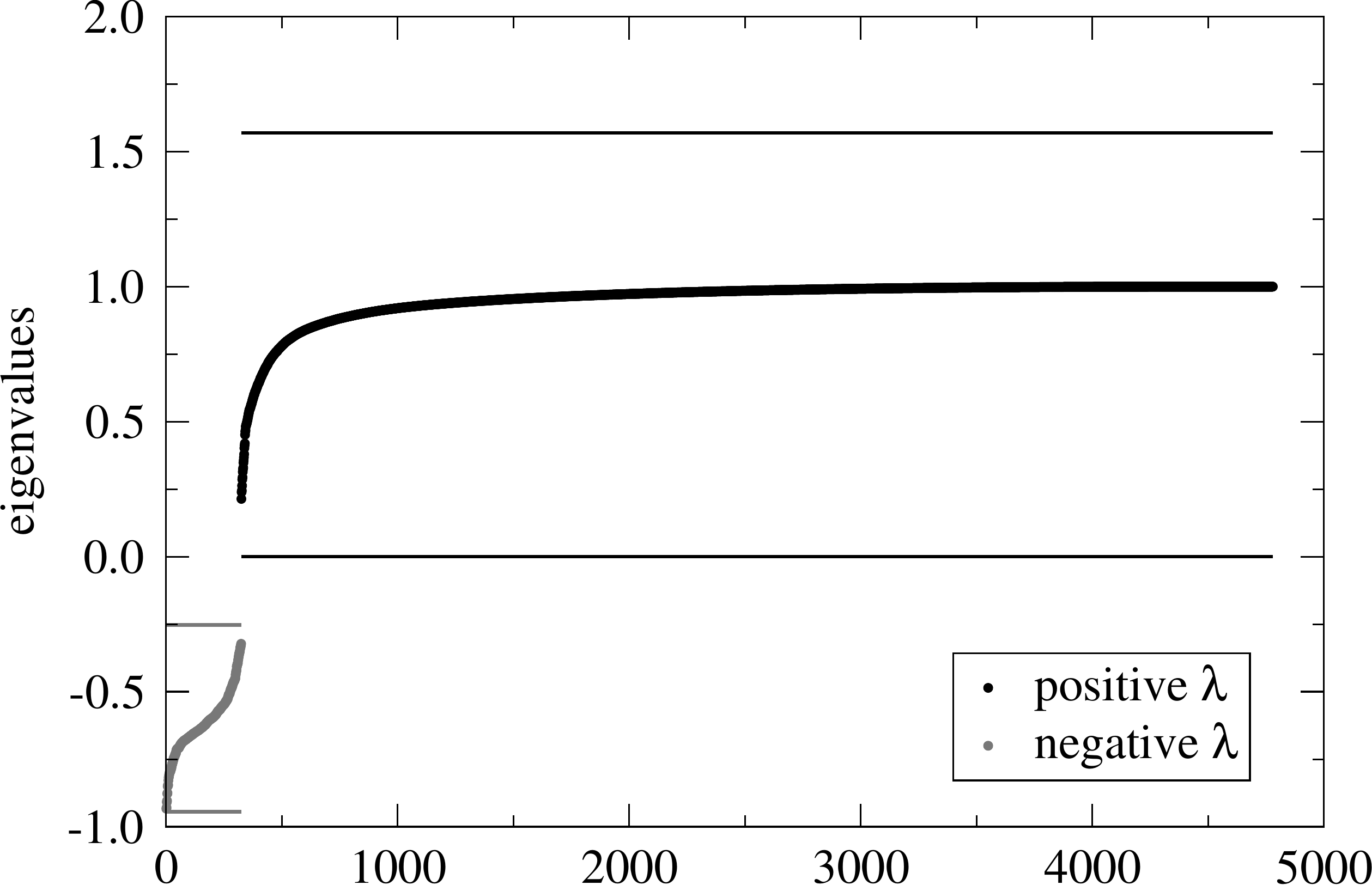}}
  \hfill\null
  \caption{Eigenspectrum of $\mathcal{M}^{-1}\mathcal{A}$ (left) and $\mathcal{M}_a^{-1}\mathcal{A}$ (right) for the test case of Figure \ref{fig:cube} 
  using $\tS_u^{-1} = \texttt{AMG}(S_u)$. In both frames the theoretical bounds of Theorem \ref{thm:sym_bounds} and \ref{thm:nonsym_bounds} are also given.} 
  \label{fig:eigsMinvA}
\end{figure}

\begin{remark}
\label{rem:loose_bounds}
Figure \ref{fig:eigsPinvA} shows that the theoretical bounds introduced in Section \ref{sec:bounds} can be indeed pretty loose in the case $\tS_u^{-1}=S_u^{-1}$. 
If we replace the application of the exact inverse of $S_u$ with an approximation, given for instance by a standard Algebraic Multigrid (AMG) algorithm,
complex eigenvalues with $\mathcal{M}^{-1}\mathcal{A}$ and positive eigenvalues with $\mathcal{M}_a^{-1}\mathcal{A}$ immediately arise (Figure \ref{fig:eigsMinvA}). 
In this case, the theoretical bounds, also provided in Figure \ref{fig:eigsMinvA}, appear to be more significant.
\end{remark}

From the observations above, we can conclude that the proposed choices for $C$ appear to be equally effective for the overall preconditioner performance. Since
equation \eqref{eq:Cdef} can be computed independently of the regularity of $A$, it should be preferred for the sake of robustness and generality. 

Another issue concerns whether it is possible to find an optimal value for the relaxing factor $\omega$.
It is clear that the selection of $\omega$ should take into account a balance between two opposite
requirements. On the one hand, we want the matrix $\hat{\mathcal{A}}$ and $\overline{\mathcal{A}}$ of equation \eqref{eq:matAhat} and \eqref{eq:matAbar}, respectively, to be as close
as possible to $\mathcal{A}$, which implies $\omega \rightarrow 0$. On the other hand, it is important that the condition number of $S_u$ does not deviate too much from that of $A$, or, if $A$ is singular, is not too large, which is in contrast with the requirement $\omega \rightarrow 0$.


\begin{figure}
  \centering
  \null\hfill
  \includegraphics[width=0.49\linewidth]{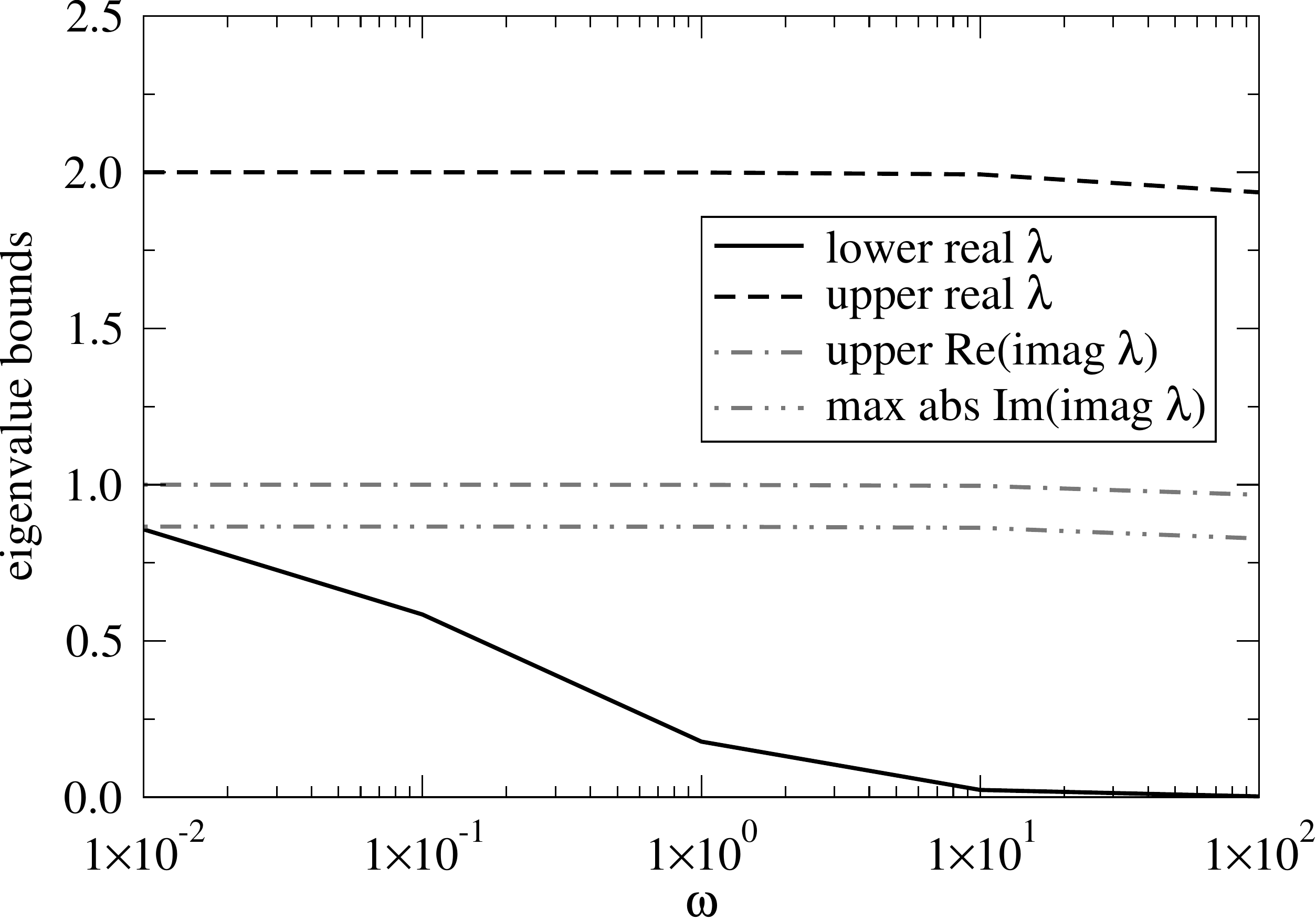}\hfill
  \includegraphics[width=0.49\linewidth]{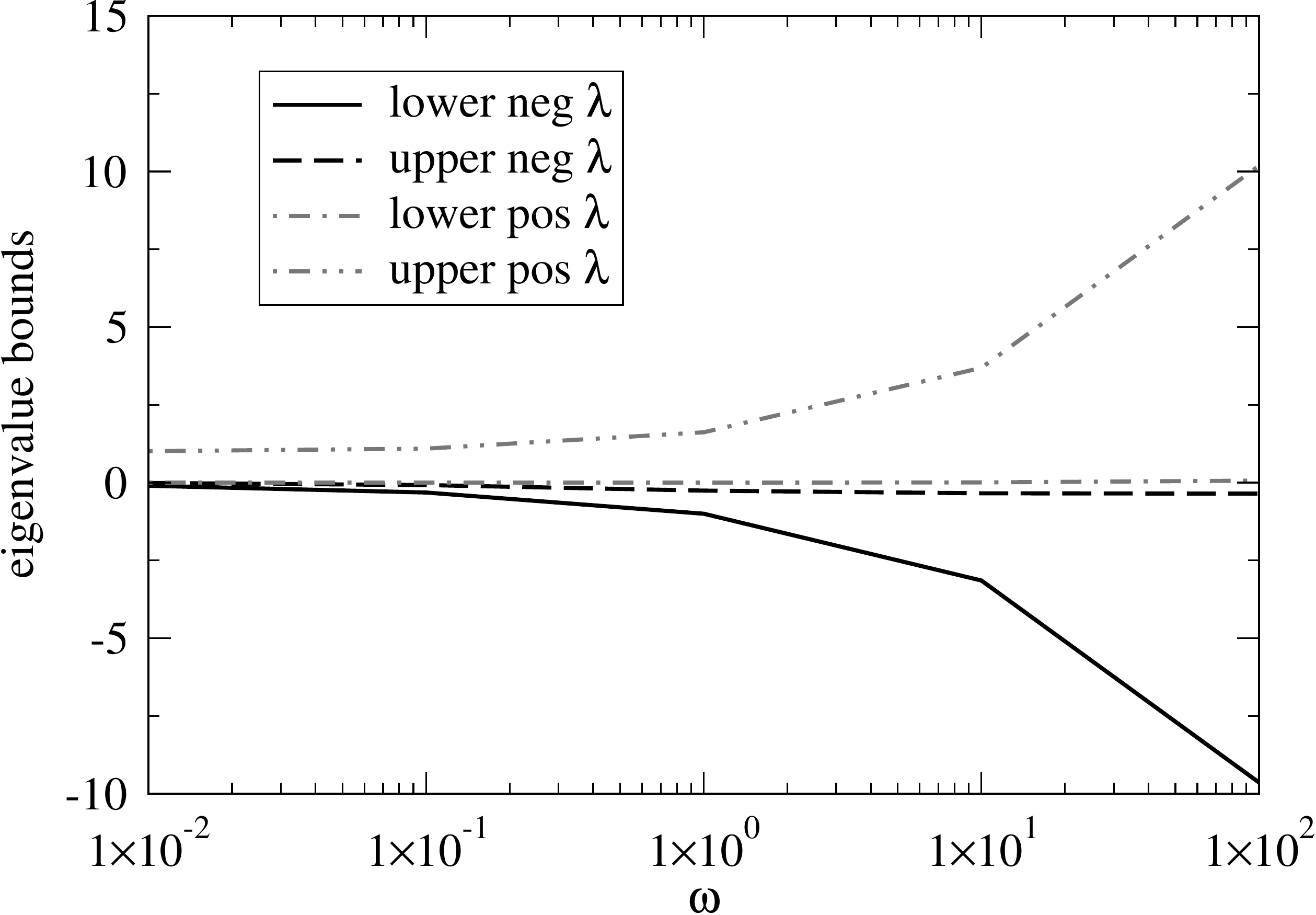}
  \hfill\null
  \caption{Theoretical eigenvalue bounds for $\mathcal{M}^{-1}\mathcal{A}$ (left) and $\mathcal{M}_a^{-1}\mathcal{A}$ (right) vs $\omega$ for the test case of Figure \ref{fig:cube}.} 
  \label{fig:omegaCfr}
\end{figure}

\begin{figure}
  \centering
  \includegraphics[width=0.49\linewidth]{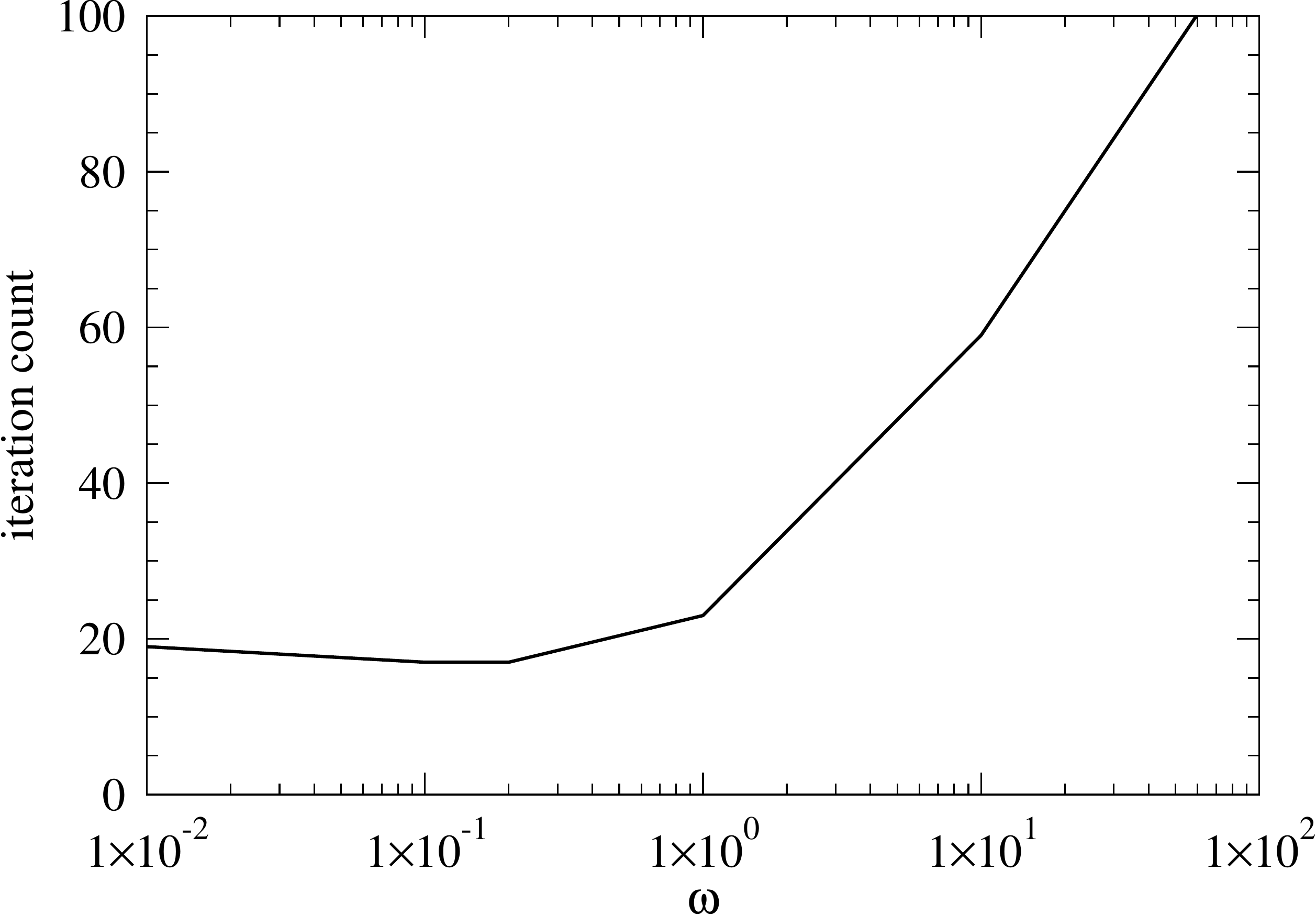}
  \includegraphics[width=0.49\linewidth]{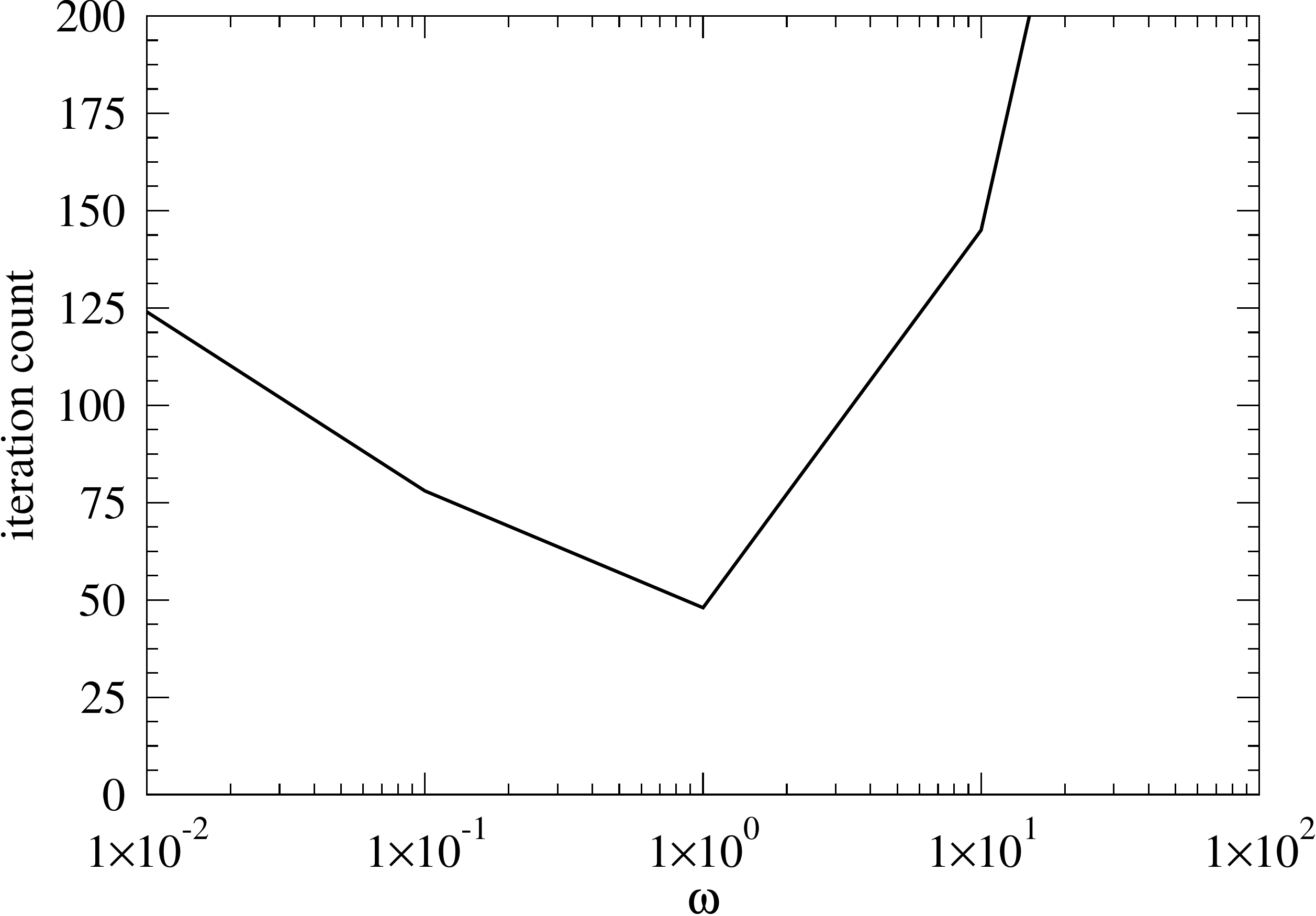}
  \caption{Iteration counts using a right-preconditioned GMRES(100) for the cube of
    Figure \ref{fig:cube} (on the left) and the \texttt{Mexico} case (on the right).}
  \label{fig:omega_iter}
\end{figure}

The role of $\omega$ is investigated numerically in the test case of Figure \ref{fig:cube}. First, we analyze the behavior of the theoretical eigenvalue bounds of Lemma \ref{lem:sym_bounds} and Theorem \ref{thm:nonsym_bounds} for $\tS_u^{-1}=S_u^{-1}$ varying $\omega$ (Figure \ref{fig:omegaCfr}).
We notice that the possible complex eigenvalues should not be affected by $\omega$ and that $\omega\leq1$ is the most convenient option, with bounds that are progressively narrower as $\omega\rightarrow0$. However, the size of the interval containing the eigenspectrum does not appear to vary significantly for $\omega\in[0.01,1]$.
In order to practically analyze the influence of $\omega$, we compute the solution of the linear system
with the matrix $\mathcal{A}$ 
preconditioned with
$\mathcal{M}$ 
for different $\omega$ values. We use a right-preconditioned GMRES(100) \cite{saad1986gmres}. Two test cases are considered: (i) the
stretched cube of Figure \ref{fig:cube}, and (ii) the \texttt{Mexico} case, i.e., a faulted reservoir discussed in the sequel
(see Section \ref{sec:numres} for a description). In Figure \ref{fig:omega_iter}, the number of
iterations is reported for different values of $\omega$ in the range $[10^{-2},10^2]$. As
expected, with $\omega > 1$ there is a clear deterioration of the
preconditioner effectiveness, while for values lower than 1 the situation is different in the
two cases. Consistently with the observations derived from Figure \ref{fig:omegaCfr}, for the
cube of Figure \ref{fig:cube} there is almost no difference in the iterations count. Instead,
in the \texttt{Mexico} case, there is an appreciable increase in the GMRES iterations when
$\omega \ne 1$. From these observations, we can conclude that, even if an optimal
$\omega$ value close to $1$ might exist, from a practical viewpoint it is not really necessary to look for it. In fact, if it exists the
iterations reduction is generally small, otherwise the loss can be significant.

To conclude this analysis, we compare the performance of the two preconditioners, $\mathcal{M}$ and $\mathcal{M}_a$,
that mainly differ in the nature of the eigenspectrum of the preconditioned matrix. In fact, as already observed, the former provides
complex eigenvalues with positive real part, while the latter yields 
real eigenvalues with both positive and negative signs. We use $C$ computed according
to Equation \eqref{eq:Cdef} with $\omega = 1$. To compare the two approaches, we use a
right-preconditioned GMRES(100) for the \texttt{Mexico} test case. Using $\mathcal{M}$,
GMRES requires 48 iterations to achieve convergence, while with $\mathcal{M}_a$, as many as 109
iterations are needed, i.e., more than twice. Hence, 
the use of $\mathcal{M}$ appears to be in practice preferable.

\section{Implementation for high performance computers}
\label{sec:impl}

An efficient parallel implementation of \precname is fundamental for solving large size problems that could arise in several applications. Most of the operations that need to be parallelized in \precname, such as approximating the (1,1) block or performing sparse matrix by matrix products, have been already studied by several authors and can be borrowed from existing libraries. In this work, we build \precname relying on the parallel Chronos software package~\cite{IsoFriSpiJan21}. Chronos is a library designed to solve linear algebra problems on high performance computers, providing the parallel preconditioners and sparse linear algebra kernels required to form and approximate the Schur complement. Chronos use is free for research purposes
and its license can be requested at the library website~\cite{CHRONOS-webpage}.

We use the distributed sparse matrix (DSMat) structure in Chronos to store the system matrix with the care of preserving the native 4-block structure of the saddle-point matrix. To this aim, we evenly partition among the processes all the rows of both the (1,1) block and the (2,2) block, with $B$ and $B^T$ following the same row subdivision. 
To reduce the communication burden in large size problems, we try to find an optimal or nearly optimal partitioning of the system matrix. Since Chronos provides a direct interface to the ParMETIS library~\cite{pMETIS}, it is relatively easy to determine a good permutation, $P_1$, for the $A$ block, such that the first chunk of rows of $A$ is assigned to process one, the second one to process two and so on. $P_1$ is also used to permute the rows of $B$ (columns of $B^T$). However, to reduce as much as possible the data exchange among processes, we need to permute the columns of $B$ (rows of $B^T$) as well. This second permutation, $P_2$, is found using a simple heuristic. By exploring the rows of $B^T$, we determine for each Lagrange multiplier the list of connected displacement degrees of freedom, hence 
the list of possible processes that could get the ownership of the corresponding row. Initializing to zero the count of Lagrange multipliers assigned to each process, we then perform a loop over all the Lagrange multipliers and assign each one 
to the process that has the lowest number of multipliers at the time of the assignment. This procedure is clearly sequential, however it can be also carried out concurrently by each process on the chunk of rows it has been temporarily assigned to. The resulting partitioning is able to reduce inter-process communications and turns out to be reasonably well-balanced at the same time. Algorithm~\ref{alg:partLag} sketches the parallel implementation of the above procedure. Note that $n_p$ denotes the number of processes and the outermost loop is distributed among processes. The symbol $\mathcal{L}_i$ denotes the set of Lagrange Multipliers assigned to the $i$-th partition, $|\cdot|$ its cardinality, and $\mbox{first}\_{i_p}$ and $\mbox{last}\_{i_p}$ are the first and last Lagrange Multipliers assigned to the initial tentative partition $i_p$, respectively.

\begin{algorithm}[t!]
\caption{Distributing Lagrange Multipliers among processes}
\begin{algorithmic}[1]
\ForAll{$i_p \in \{ 1,\dots,n_p \}$}
   \For {$i = 1,\dots,n_p$}
      \State {$\mathcal{L}_i=\{ \emptyset \}$;}
   \EndFor
   \For{$l = \mbox{first}\_{i_p},\dots,\mbox{last}\_{i_p}$}
      \For{$j = 1,\dots,n$}
        \If {$B_{ji} \ne 0$}
            \If {$\exists \; k : j \in \mathcal{L}_k$}
               \If{$|\mathcal{L}_j| < |\mathcal{L}_k|$}
                  \State{$\mathcal{L}_j = \mathcal{L}_j \cup \{ l\}$};
               \EndIf
            \Else
               \State{$\mathcal{L}_j = \mathcal{L}_j \cup \{ l\}$};
            \EndIf
         \EndIf
      \EndFor
   \EndFor
\EndFor
\end{algorithmic}
\label{alg:partLag}
\end{algorithm}

\begin{figure}
  \centering
  \includegraphics[height=0.4\linewidth]{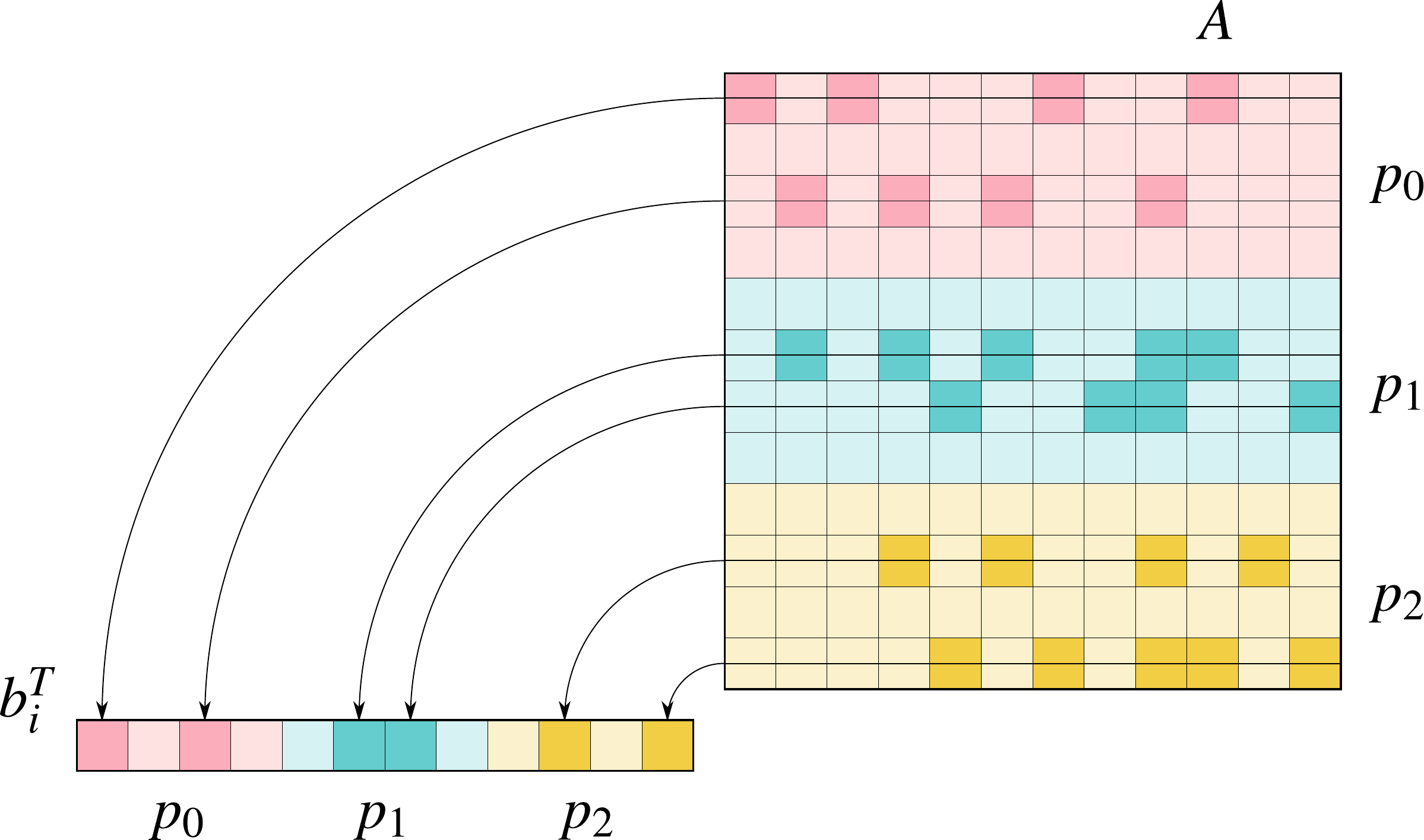}
  \caption{Schematic representation of the data exchange required for the computation of the $i$-th diagonal entry of $C$. The matrix $A$ is partitioned between three processes, $p_0$, $p_1$ and $p_2$, and the red lines highlight the rows of $A$ that must be sent to the process owning the $i$-th Lagrange multiplier. Non-zero entries in $\mathbf{b}_i^T$ show the row indices that need to be collected from $A$.}
  \label{fig:C_comp}
\end{figure}

Once the problem is adequately partitioned, some communication is needed to compute the diagonal matrix $C$. Indeed, each entry of $C$ is computed according to~(\ref{eq:Cdef}) with a few rows of $A$ exchanged among processes to form the submatrix $A|_{b_i}$. This data exchange is very similar to the one required in a matrix-by-matrix product between $B^T$ and $A$. In fact, assuming that process $p$ is the owner of the $i$-th Lagrange multiplier, according to Figure~\ref{fig:C_comp}, it must receive from the other processes the rows of $A$ corresponding to the non-zero columns in $\mathbf{b}_i^T$, the $i$-th row of $B^T$. After the collection of all the $A$ rows required by the corresponding Lagrange multipliers, process $p$ can locally form each submatrix $A|_{b_i}$ and compute the entries $C_{i,i}$.

The remaining part of the \precname parallelization only requires direct calls to Chronos functions. In fact, the Schur complement is simply formed by calling the distributed sparse matrix-by-sparse matrix product, which manages all the necessary communications on its own. As a preconditioner for $S_u$, we choose AMG with the special set-up designed for structural problems~\cite{IsoFriSpiJan21}, which includes a Factorized Sparse Approximate Inverse (FSAI) as smoother and adaptive 
interpolation.
The source files needed to store the saddle-point matrix, compute the \precname preconditioner and solve the linear system can be found at \url{https://github.com/matteofrigo5/HPC\_ReverseAugmentedConstrained} and they can be easily modified to test other approaches or use another preconditioner for $S_u$.

\section{Numerical results}
\label{sec:numres}

The purpose of this section is to investigate the numerical properties
of the proposed approach in relatively simple, though challenging, test cases, and the related computational performance and parallel efficiency in real-world applications. First of
all, four medium-size examples are analyzed, in order to understand the properties of the
different matrices involved in the preconditioner definition, the influence of the ratio
between $n_t$ and $n_u$, and the sensitivity to different discretizations and contact modeling approaches.
Then, two large-size real-world cases are presented to shed light on the \precname computational
effectiveness.

\subsection{Theoretical analysis}

\begin{table}
  \centering
  \begin{tabular}{l|r|r|r|l}
    case & $n_u$ & $n_t$ & $n_t/n_u$ & comments \\
    \hline
    \texttt{floating-side} & 218,790 &   6,336 &  2.9\% & \textit{fault-constrained}
      simple structured case \\
    \texttt{Mexico}        & 171,150 &  12,027 &  7.0\% & real world subsidence model
      \cite{franceschini2015modelling} \\
    \texttt{15flt}         & 379,983 & 167,799 & 44.2\% & large 2D/3D ratio
      \cite{franceschini2019block} \\
    \texttt{node-surf}     & 194,208 &   6,936 &  3.6\% & node-to-surface contact formulation \\
  \end{tabular}
  \caption{Matrix sizes for the analyzed test cases.}
  \label{tab:sizeCases}
\end{table}

\begin{figure}
  \centering
  \null\hfill
  \vcenteredhbox{\subfloat[\texttt{floating-side}.]
    {\includegraphics[height=0.25\linewidth,trim=-100 0 -100 0]{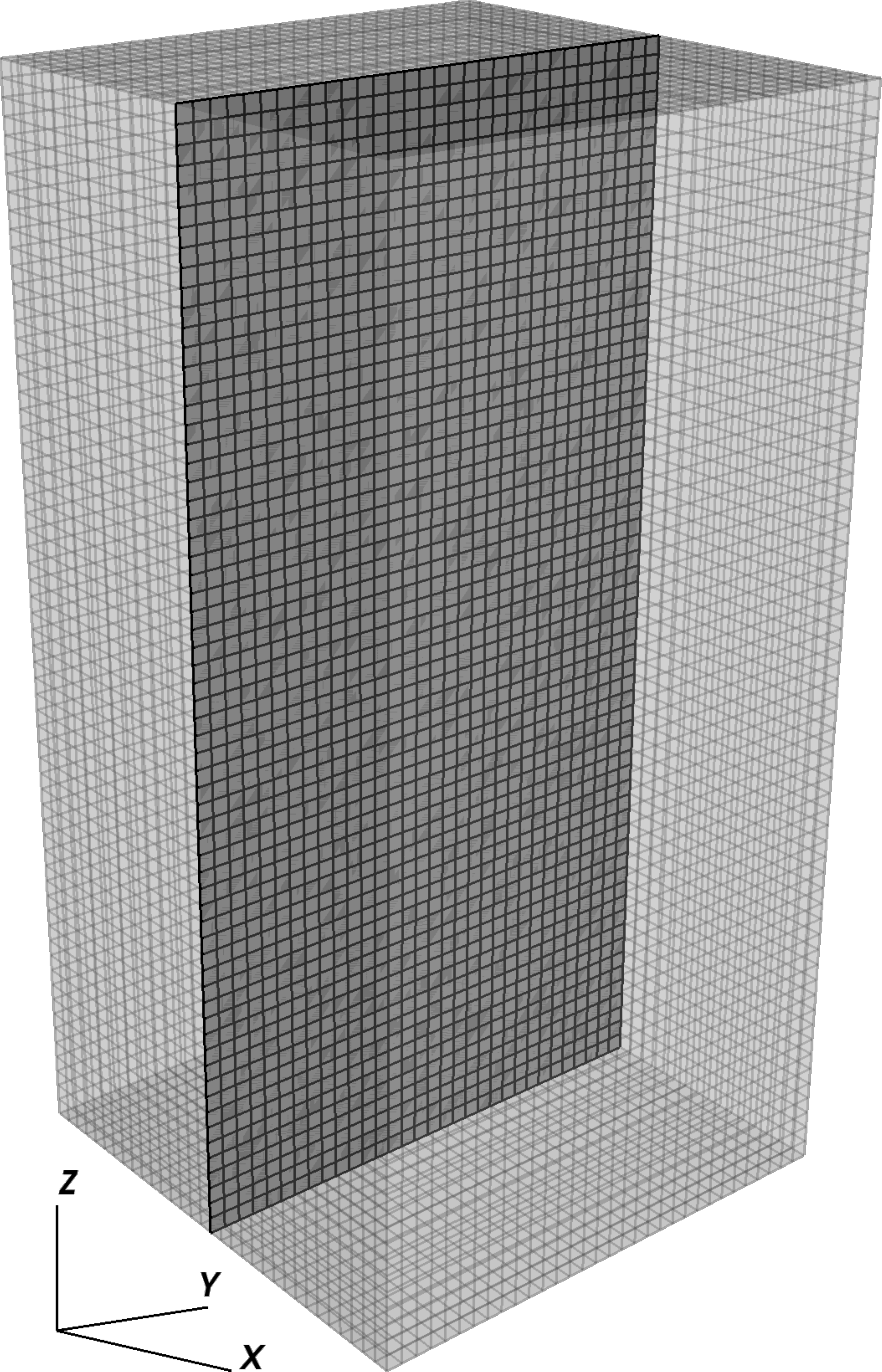}}}\hfill
  \vcenteredhbox{\subfloat[\texttt{Mexico}. The vertical exaggeration factor is 10.]
    {\includegraphics[height=0.20\linewidth]{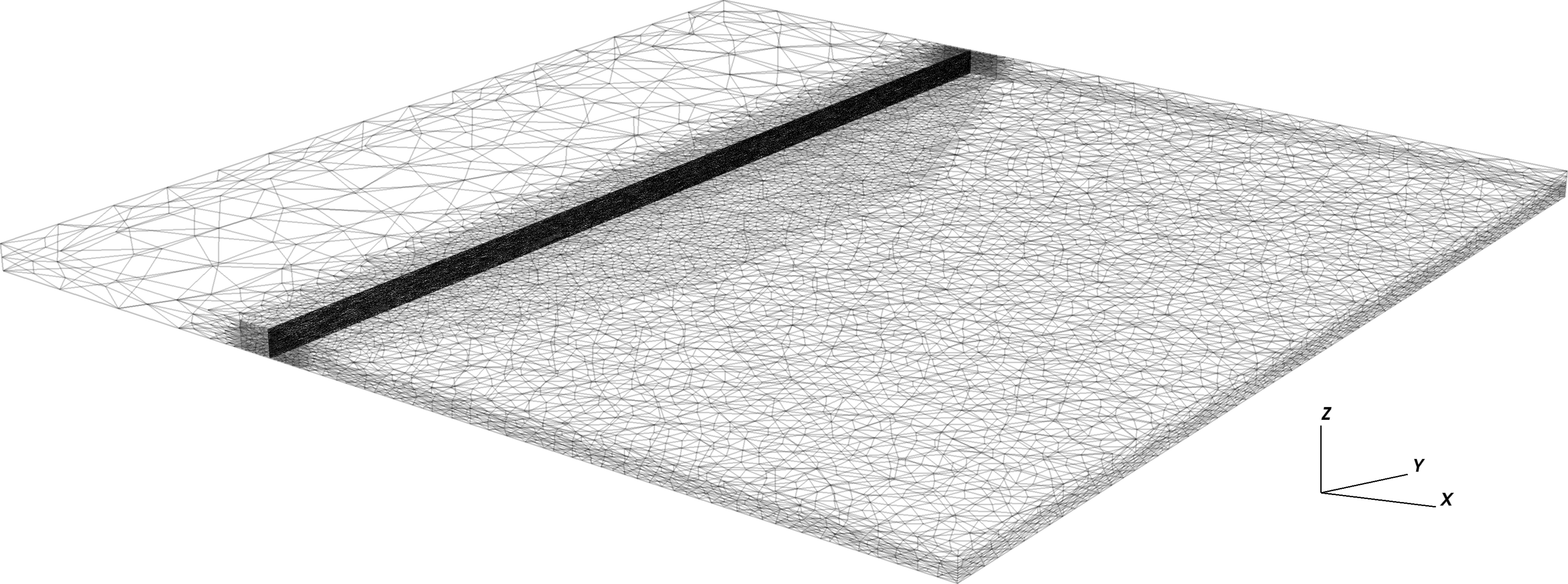}}}
  \hfill\null\\
  \null\hfill
  \subfloat[\texttt{15flt}.]
    {\includegraphics[height=0.25\linewidth]{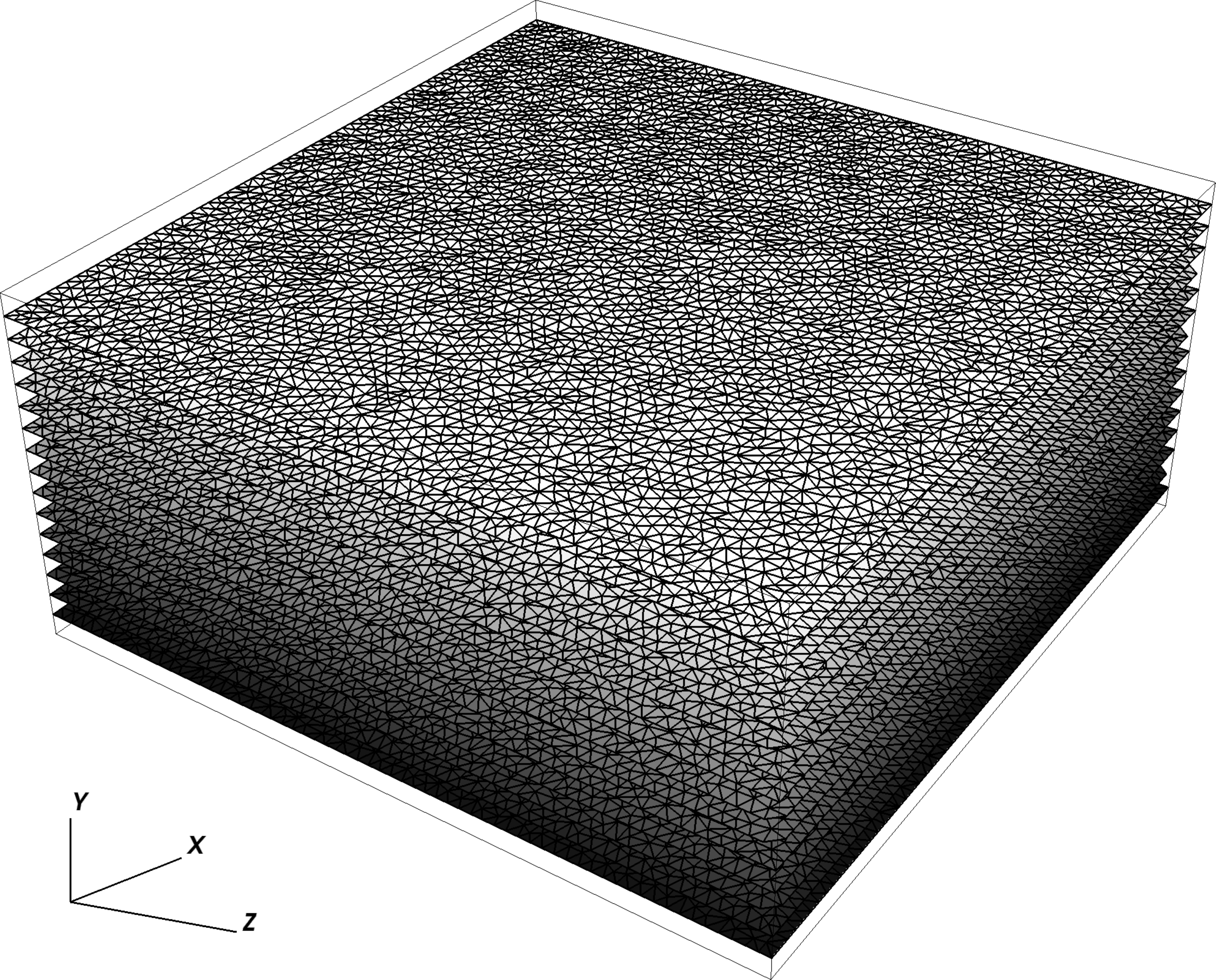}}\hfill
  \subfloat[\texttt{node-surf}.]
    {\includegraphics[height=0.25\linewidth]{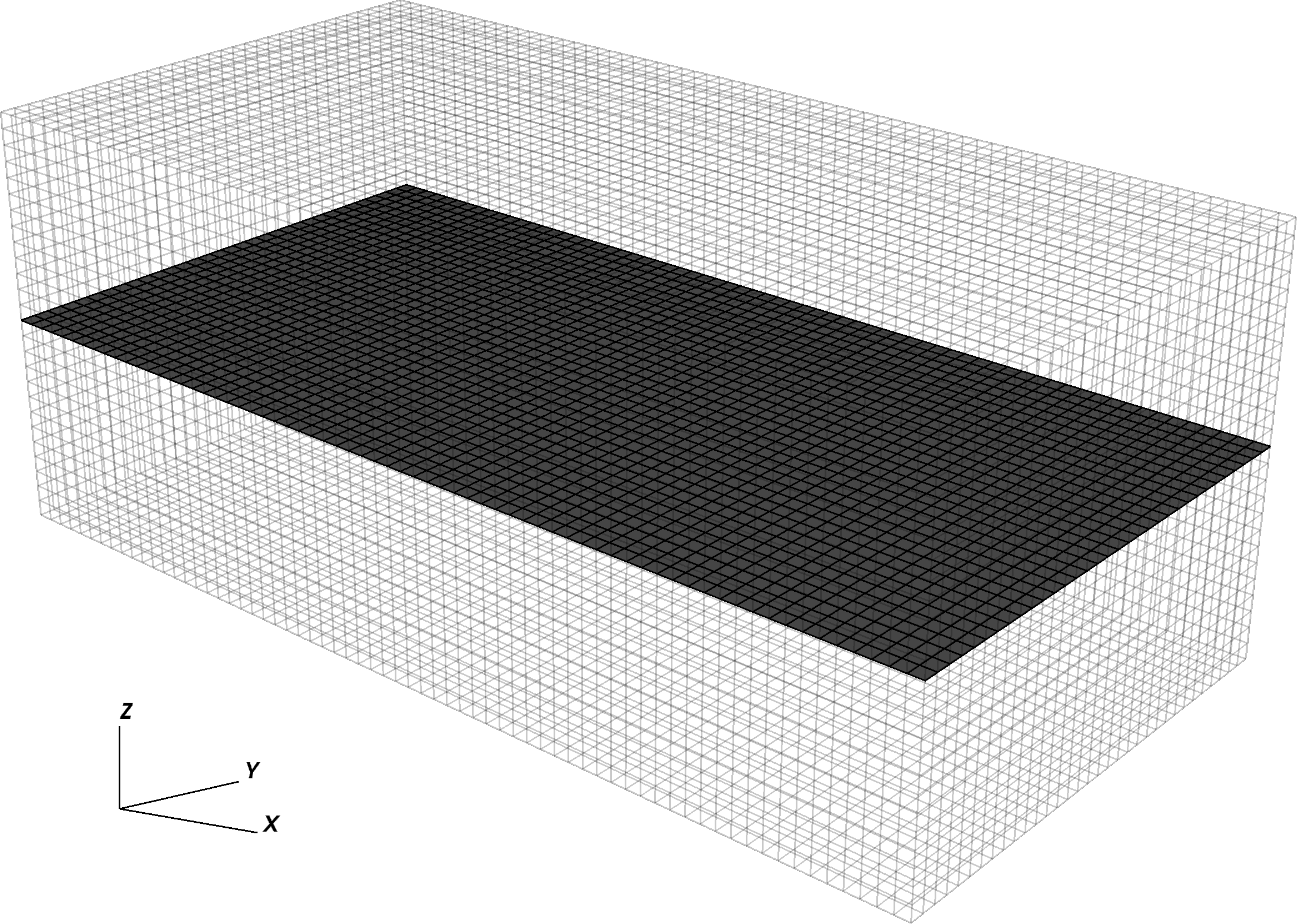}}
  \hfill\null
  \caption{Computational meshes for the test cases listed in Table \ref{tab:sizeCases}.
    Note that for the sake of visibility, only the fractures are shown for the
    \texttt{15flt} case.}
  \label{fig:meshCases}
\end{figure}

Table~\ref{tab:sizeCases} provides the size and a short description of the four matrices used
in this analysis. Figure~\ref{fig:meshCases} shows a sketch of the computational meshes with the
contact surfaces highlighted in darker gray. 
The features of the investigated problems are as follows:
\begin{enumerate}
\item {\tt floating-side} is designed so as to have a singular leading block. It has a passing fracture and
only the leftmost portion of the domain is constrained by Dirichlet conditions. The remaining part of the body is constrained by contact conditions, thus leading to a
reducible $A$ block with a singular sub-block; 
\item {\tt Mexico} originates from a model used to simulate the land
subsidence process above an over-exploited faulted aquifer system \cite{franceschini2015modelling}. The problem is characterized by a strong element distortion, with a large width in the horizontal plane and a relatively limited extension in the vertical
direction, including a vertical fracture that encompasses the whole domain; 
\item {\tt 15flt} discretizes 
a cube with 15 parallel fractures. This example is taken from the analysis developed in \cite{franceschini2019block} and is specifically designed to maximize the ratio between
$n_t$ and $n_u$, i.e., the number of traction and displacement unknowns;
\item {\tt node-surf} simulates the motion of a solid block sliding over a constrained deformable basement. 
\end{enumerate}
The first three problems have
been obtained with the aid of a \textit{node-to-node} Lagrange multipliers approach \cite{fraferjantea16}, while,
as it is also suggested by
the name, the fourth one is based on a \textit{node-to-surface} contact model.

\begin{figure}
  \centering
  \null\hfill
  \subfloat[\texttt{floating-side} case.]
    {\includegraphics[width=0.40\linewidth]{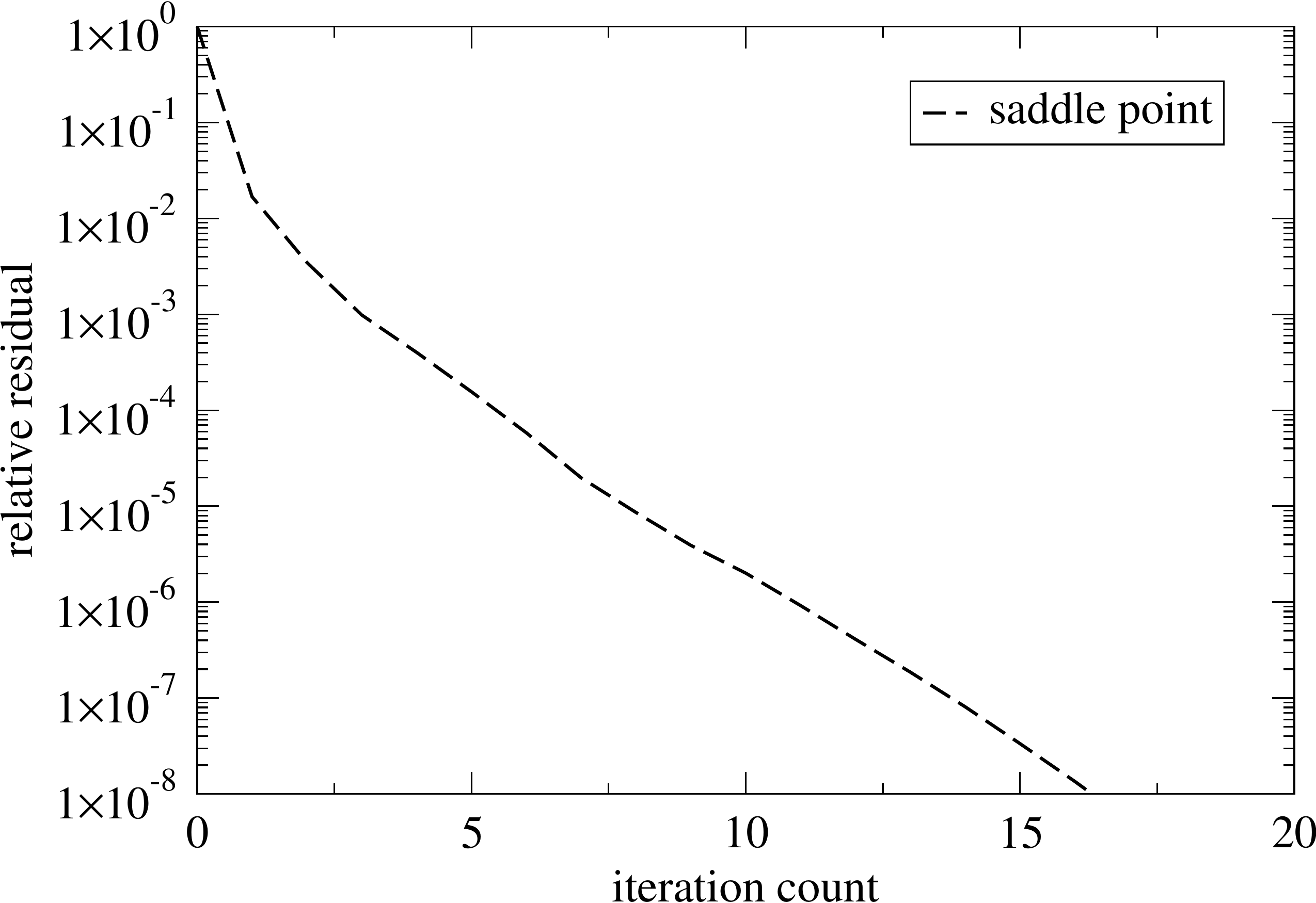}}\hfill
  \subfloat[\texttt{Mexico} case.]
    {\includegraphics[width=0.40\linewidth]{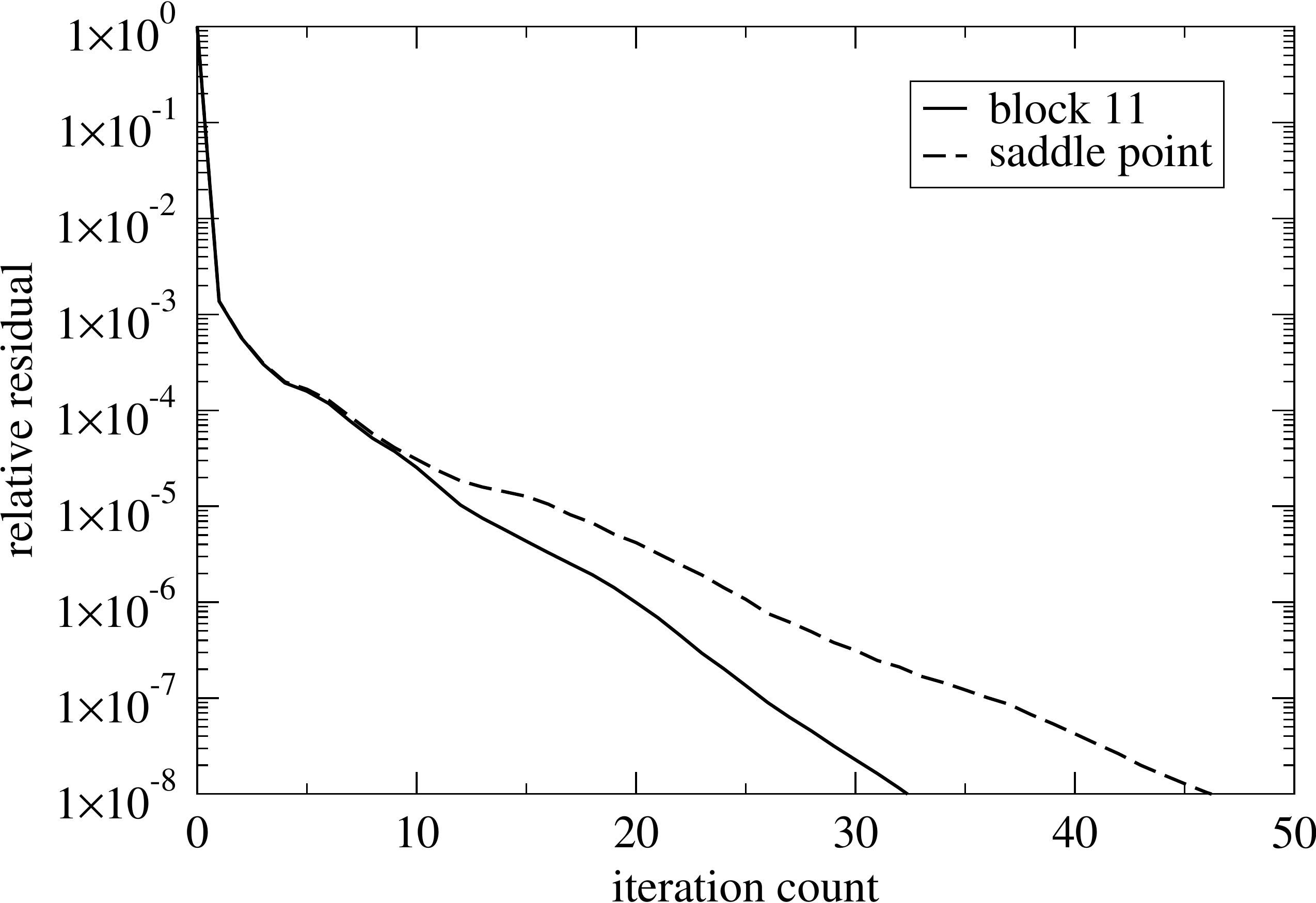}}
  \hfill\null\\
  \null\hfill
  \subfloat[\texttt{15flt} case.]
    {\includegraphics[width=0.40\linewidth]{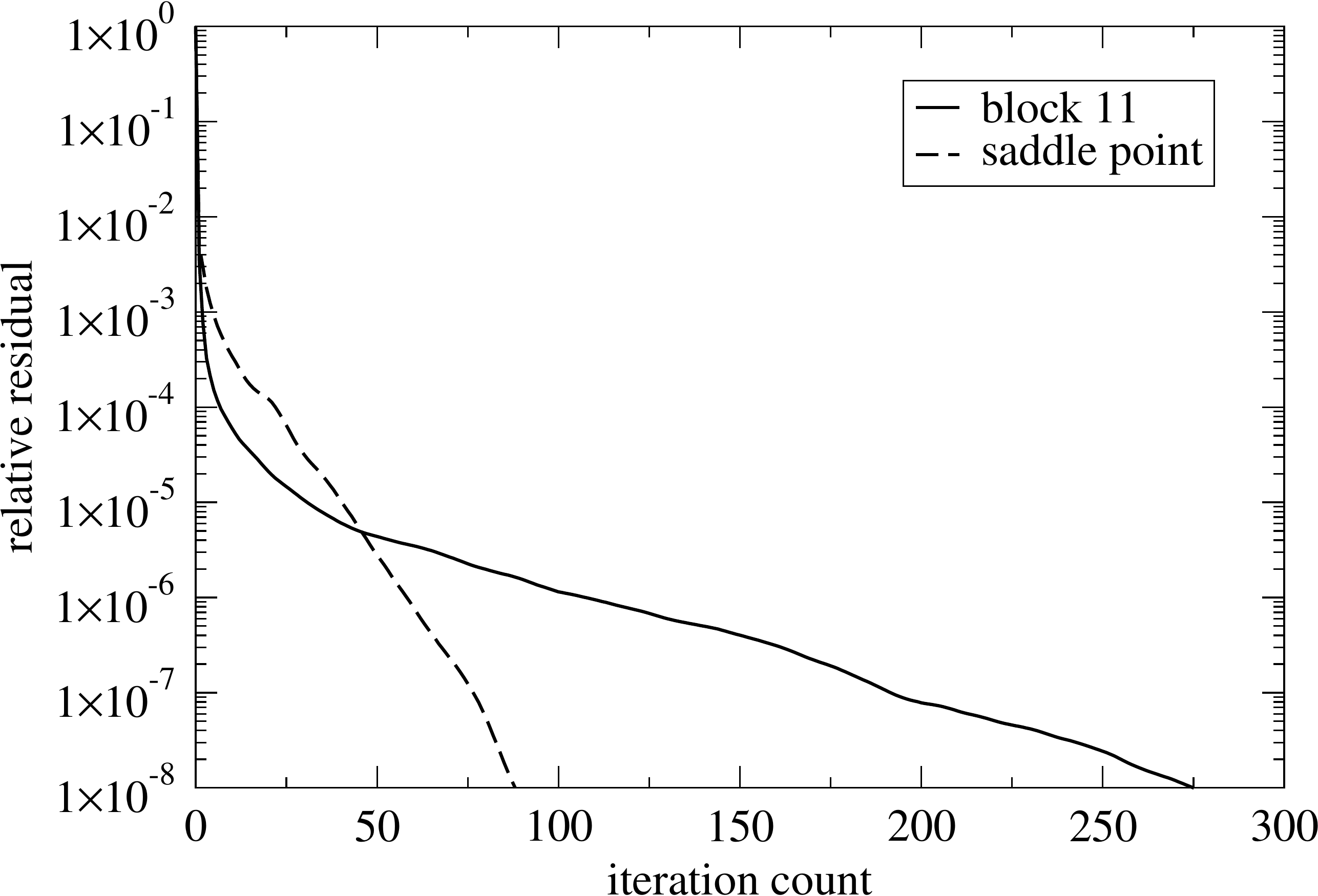}}\hfill
  \subfloat[\texttt{node-surf} case.]
    {\includegraphics[width=0.40\linewidth]{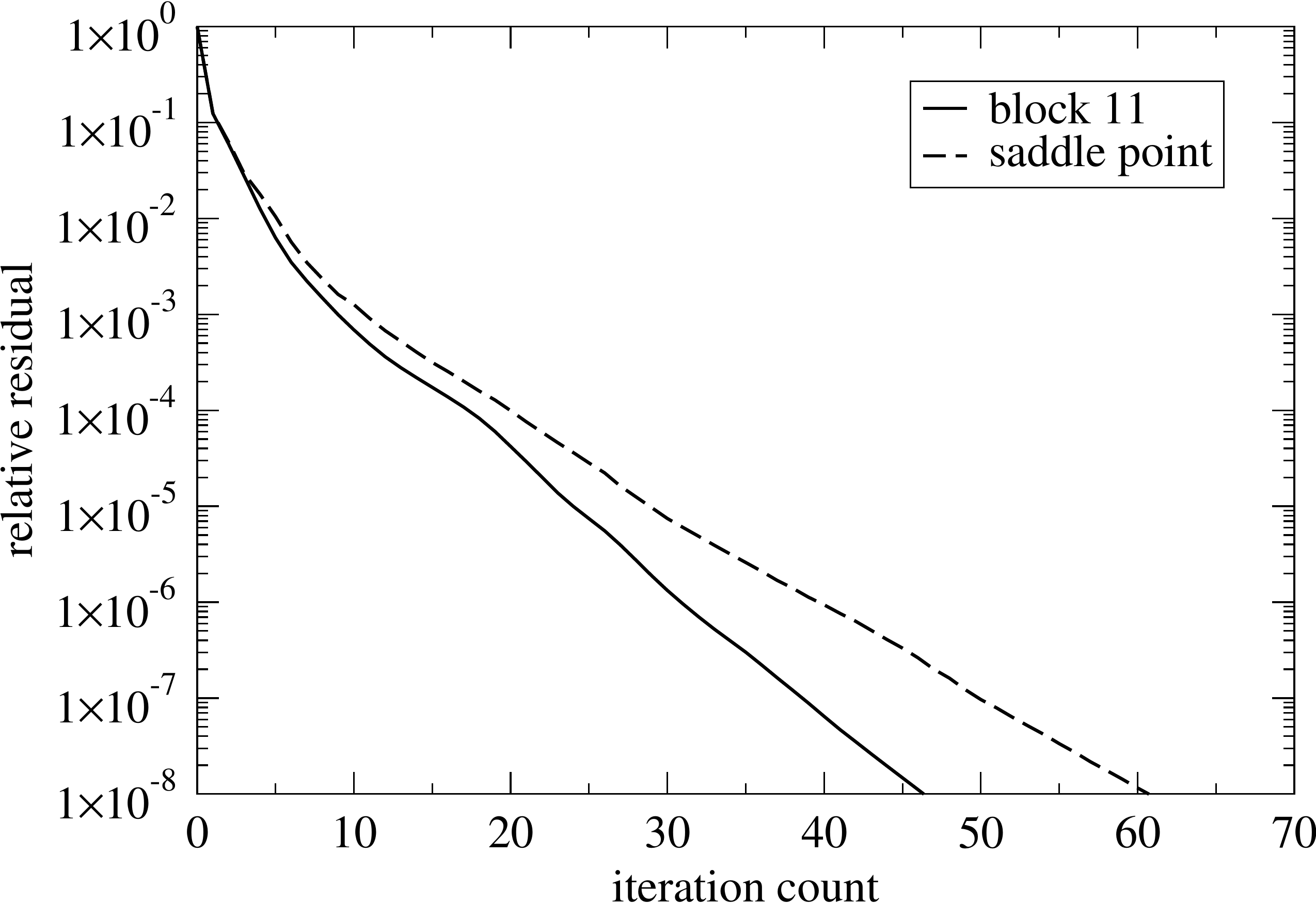}}
  \hfill\null
  \caption{Convergence profiles for both the complete saddle point matrix and the leading
    block only for all the test cases.}
  \label{fig:convCases}
\end{figure}

First, we want to compare the convergence properties in the solution of the full saddle-point system with that of the
leading block $A$ alone.
This is done for all the test cases except for the {\tt floating-side} problem, where the leading block is singular by construction and only the full saddle-point system is solved. 
A right-preconditioned GMRES(100)
is used with a unitary right-hand side. The iterations are stopped when the 2-norm of the initial residual is reduced by 8 orders of magnitude.
The \precnameVar{1} approach is used as a preconditioner, with the application of the inverse of the primal Schur complement, $\widetilde{S}_u^{-1}$, carried out approximately by the AMG method proposed in \cite{IsoFriSpiJan21,paludetto2019novel}.
For the sake of the comparison, we solve the system 
with the leading block $A$ alone by accelerating GMRES(100) with the same multigrid method 
used to approximate 
$\widetilde{S}_u^{-1}$.
Note that
this implies almost identical operator and grid complexities. 
Figure
\ref{fig:convCases} shows the convergence profiles for all the cases. 
Often, Augmented Lagrangian approaches may face the issue of worsening the conditioning of the augmented block. This is not the case, however, for the selected strategy.
While for the
\texttt{Mexico} and \texttt{node-surf} cases the convergence profiles are similar with a slight increase of the iteration counts in the saddle-point problem, they are
quite different in the \texttt{15flt} case, but in the opposite direction. In this case, the leading block alone 
turns out to be more ill-conditioned than the augmented one. 
This is what usually happens when $A$
is much less
constrained than the full saddle-point matrix.
Hence, we can conclude that the primal Schur complement obtained by augmenting the leading block has a conditioning that is similar to, or even better than, the leading block itself, and especially so if $A$ is close to be singular.
In the
\texttt{floating-side} case, we cannot perform this comparison, but we can observe that the convergence profile 
turns out to be
similar to the other test cases.

The extreme eigenvalues, $\lambda_{\min}$ and $\lambda_{\max}$, and the condition number, $\kappa$, for both the
leading block $A$ and the primal Schur complement $S_u$ are reported in Table \ref{tab:ExtEigs} for the test cases considered above. 
In Table \ref{tab:ExtEigs}, $\lambda_{\min,0}$ is the first
nonzero eigenvalue and $\kappa_{\text{eff}}$ is the effective condition number, i.e.,
the ratio between the largest and the smallest nonzero eigenvalue, that is used for problem
1 where the leading block is singular.
As expected, it can be seen that $\lambda_{\max}(\mat{A})$ and $\lambda_{\max}(\mat{S_u})$,
i.e. the maximum eigenvalues, are quite close each other.
Indeed, matrix $C$ was built so as to
be spectrally equivalent to $B^TA^{-1}B$, i.e., the eigenspectrum of the primal Schur complement $S_u$ should be aligned with the non-zero part of that of $A$.
On the other hand,
the minimum eigenvalue, $\lambda_{\min}$, generally increases from $A$ to $S_u$. In fact,
from a physical viewpoint, the contribution $B C^{-1}
B^T$ can be regarded as a sort of \textit{fracture stiffness} that is similar to that of the surrounding medium and is added to the original stiffness
matrix, thus increasing the overall problem stiffness.
As a consequence, $S_u$ has often a conditioning similar to $A$, or is even better conditioned than $A$. 
Note that the largest reduction in
the condition number corresponds to the largest difference in the iteration count (see
Figure \ref{fig:convCases}).



\begin{table}
  \centering
  \begin{tabular}{l|c|cc|c}
    case & matrix & $\lambda_{\min}$ $(\lambda_{\min,0})$ & $\lambda_{\max}$ &
      $\kappa$ $(\kappa_{\text{eff}})$ \\
    \hline
    \hline
    \multirow{2}{*}{\texttt{floating-side}} &
      $A$    & $(5.989\cdot 10^{-5})$ & $7.449 \cdot 10^{-1}$ & $(1.244 \cdot 10^{4})$ \\
    & $S_u$    & $3.946 \cdot 10^{-5}$ & $7.449 \cdot 10^{-1}$ & $1.888 \cdot 10^{4}$ \\
    \hline
    \multirow{2}{*}{\texttt{Mexico}} &
      $A$    & $1.347 \cdot 10^{9}$ & $5.651 \cdot 10^{14}$ & $4.197 \cdot 10^{5}$ \\
    & $S_u$    & $2.876 \cdot 10^{9}$ & $5.651 \cdot 10^{14}$ & $1.965 \cdot 10^{5}$ \\
    \hline
    \multirow{2}{*}{\texttt{15flt}} &
      $A$    & $5.173 \cdot 10^{1}$ & $2.870 \cdot 10^{8}$ & $5.549 \cdot 10^{6}$ \\
    & $S_u$    & $3.530 \cdot 10^{3}$ & $3.603 \cdot 10^{8}$ & $1.021 \cdot 10^{5}$ \\
    \hline
    \multirow{2}{*}{\texttt{node-surf}} &
      $A$    & $6.439 \cdot 10^{0}$ & $1.364 \cdot 10^{6}$ & $2.119 \cdot 10^{5}$ \\
    & $S_u$    & $1.286 \cdot 10^{1}$ & $1.365 \cdot 10^{6}$ & $1.061 \cdot 10^{5}$ \\
  \end{tabular}
  \caption{Limiting eigenvalues for the leading block and the Schur complement $S_u$ for the test cases of Table \ref{tab:sizeCases}.}
  \label{tab:ExtEigs}
\end{table}

The iteration count, $n_{it}$, with the preconditioner application cost, $c_{\text{app}}$, is reported for each test case in Table
\ref{tab:cfrMCP}. The application cost denotes the number of
floating point operations required to apply the preconditioner with respect to the
matrix-vector product with the matrix $\mathcal{A}$. 
The total solution cost, $C_s$, is then estimated as the overall number of equivalent matrix-by-vector products with $\mathcal{A}$ needed to achieve the desired accuracy.
As a comparison, we consider the performance obtained
with other methods available from the literature. Table
\ref{tab:cfrMCP} reports the iteration counts, the preconditioner application costs and the total solution costs for a Mixed
Constrained Preconditioner (MCP) \cite{bergamaschi2008mixed,franceschini2019block} where:
\begin{itemize}
    \item the inverse of the leading block $A$ is approximately applied by using either the same AMG approach used with \precnameVar{1} \cite{IsoFriSpiJan21,paludetto2019novel}, or an adaptive FSAI (aFSAI) \cite{janna2015fsaipack,franceschini2019block};
    \item the Schur complement $S=-B^TA^{-1}B$ is computed inexactly by replacing $A^{-1}$ with an aFSAI of $A$ and the application of its inverse is carried out approximately by another aFSAI. 
\end{itemize}
Depending on the choice for the inner preconditioner of the leading block, we distinguish between MCP+AMG
and MCP+FSAI in Table \ref{tab:cfrMCP}. 
For these approaches, the set of user-specified
parameters providing the empirically optimal performance, as suggested in the relevant literature \cite{janna2015fsaipack, paludetto2019novel,
franceschini2019robust,IsoFriSpiJan21}, is used. It can be noticed that \precname is able to solve all the test
cases with acceptable computational costs. By distinction, a standard MCP approach cannot be used for the {\tt floating-side} test case, which has a singular leading
block. Moreover, the combination MCP+AMG fails in achieving convergence after 1,000 iterations also in the {\tt 15flt} test case, while
MCP+FSAI cannot solve the {\tt node-surf} test case. 
Even when MCP works,
the cost is usually larger than \precnameVar{1}. The reason stems from both the high AMG cost, which is required twice in the MCP+AMG approach, and the generally larger iteration count. 
Only in the \texttt{Mexico}
case the MCP+FSAI method outperforms \precname, mainly because the problem is small enough to be
tackled much more efficiently by a single level preconditioner like aFSAI.
Hence, we can conclude that the \precnameVar{1} approach appears to be much more robust and generally more efficient than MCP.

\begin{table}
  \centering
  \begin{tabular}{l|rrr|rrr|rrr}
    \multirow{2}{*}{case} & \multicolumn{3}{c|}{\precnameVar{1}} & \multicolumn{3}{c|}{MCP+AMG} &
      \multicolumn{3}{c}{MCP+FSAI} \\
    \cline{2-10}
    & $n_{it}$ & $c_{\text{app}}$ & $C_s$ & $n_{it}$ & $c_{\text{app}}$ & $C_s$ &
      $n_{it}$ & $c_{\text{app}}$ & $C_s$ \\
    \hline
    \texttt{floating-side} & 17 & 5.36 & 108.12 & --- &  --- &     --- & --- &  --- &     --- \\
    \texttt{Mexico}        & 47 & 4.37 & 252.39 & 174 & 8.76 & 1697.51 &  78 & 1.67 &  208.26 \\
    \texttt{15flt}         & 89 & 4.52 & 491.28 &   * & 9.51 &     --- & 270 & 3.18 & 1128.60 \\
    \texttt{node-surf}     & 61 & 4.21 & 317.81 & 153 & 8.43 & 1442.79 &   * & 1.85 &     --- \\
  \end{tabular}
  \caption{Computational costs and comparison with other solvers. $*$ means that full GMRES does not reach
    convergence within 1,000 iterations. --- means that the value cannot be computed.}
  \label{tab:cfrMCP}
\end{table}

\subsection{Real-world applications}

The effectiveness and the computational performance of the proposed preconditioner is finally tested in two challenging large-size applications. 

\begin{figure}
  \centering
  \null\hfill
  \includegraphics[height=0.23\textwidth]{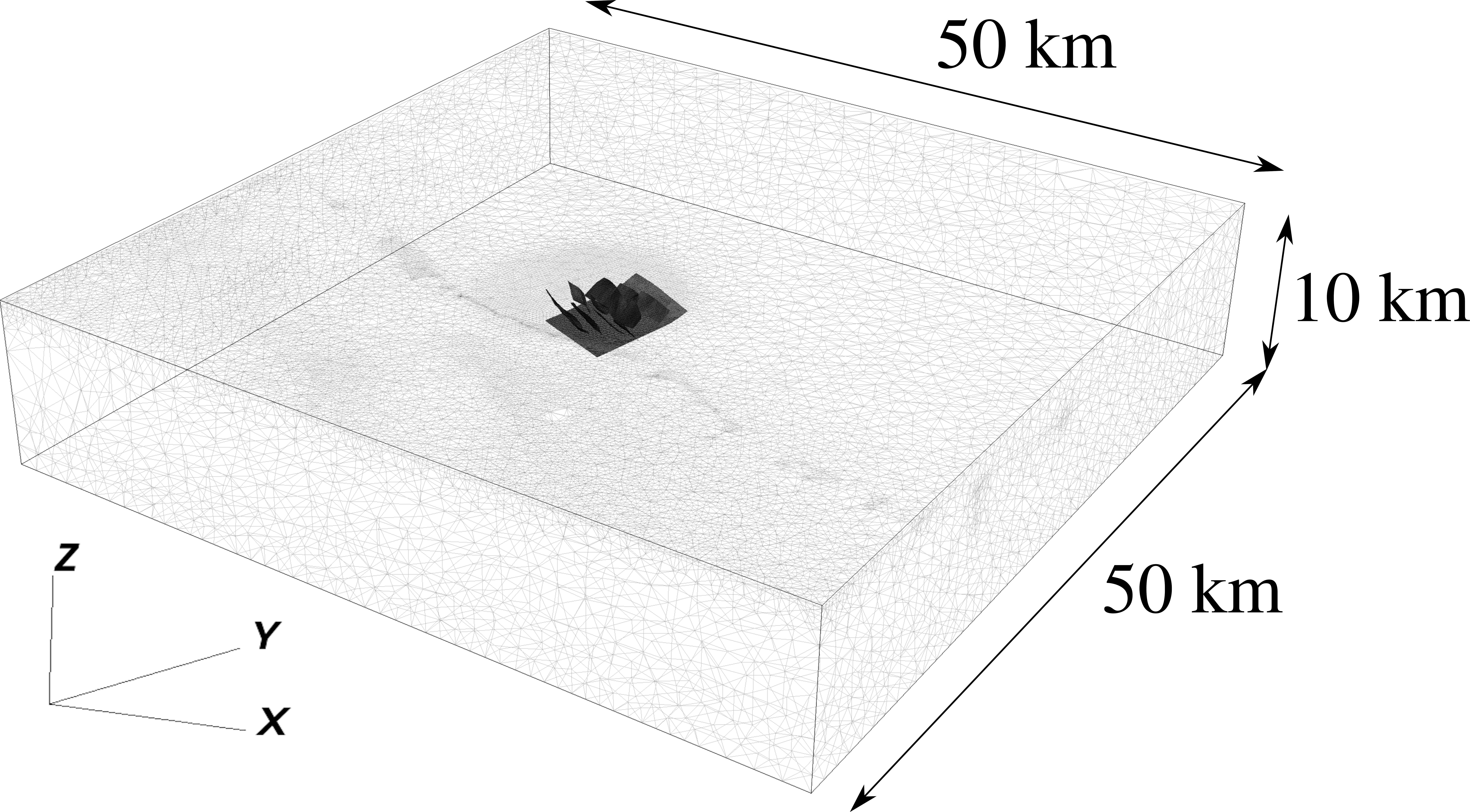}\hfill
  \includegraphics[height=0.23\textwidth]{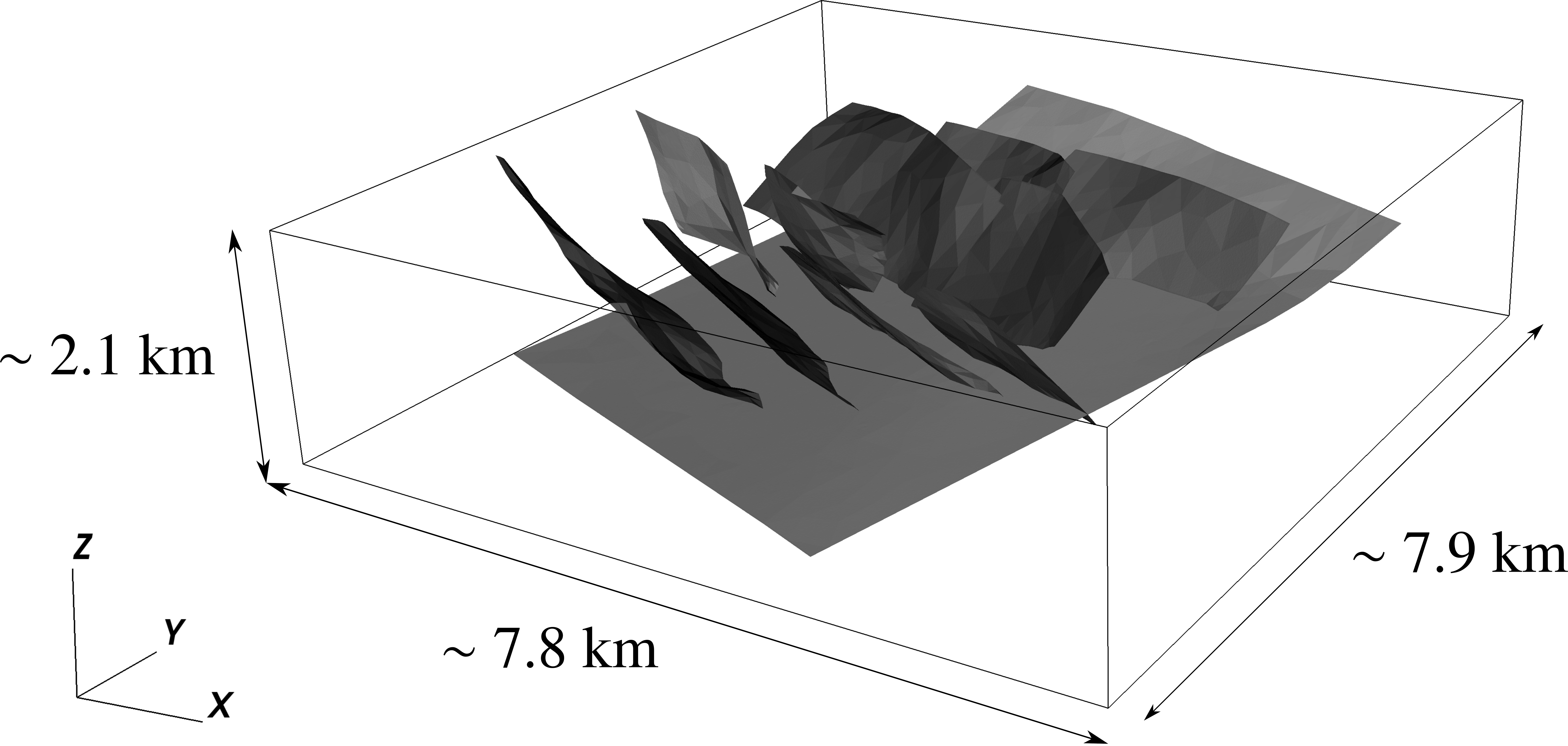}
  \hfill\null
  \caption{Test A: Computational mesh (left) with a zoom with the faults around the reservoir (right).}
  \label{fig:AmedeaMesh}
\end{figure}

\begin{figure}
  \centering
  \includegraphics[width=0.50\textwidth]{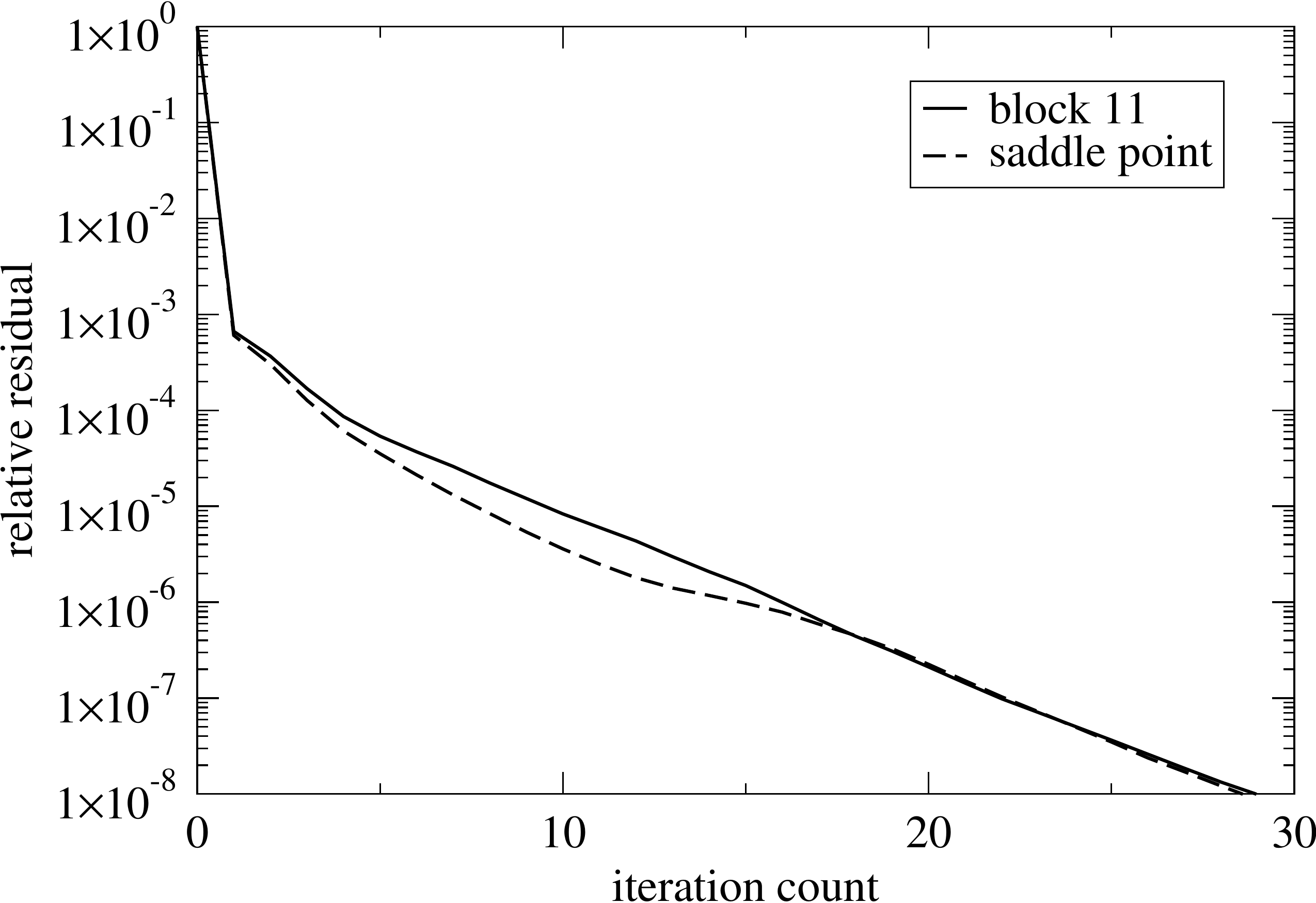}
  \caption{Test A: Convergence profile for leading block and saddle point matrix.}
  \label{fig:AmedeaConv}
\end{figure}

The first application, denoted as {\tt Test A}, arises from a real-world geomechanical model of CO$_2$ geological sequestration in a deep faulted sedimentary formation located in Italy \cite{castelletto2013geological, castelletto2013multiphysics}. This very challenging problem combines the grid distortion due to the complex geometry with the presence of 11 variably oriented faults located close to the injection wells (Figure \ref{fig:AmedeaMesh}). The purpose of the simulation was assessing the long-term safety of the sequestration with respect to possible leakages induced by the fault mechanical activation.
The mesh has 20,061,710 tetrahedral elements, 1,237,077
interface elements along the faults and 3,842,853 nodes. The problem sizes are $n_u = 11,528,559$, $n_t =
1,874,436$, with a number of non-zeros of the matrices equal to $\text{nnz}(A) = 499,168,107$ and $\text{nnz}(B) = 11,246,616$. The ratio $n_t/n_u$
is 16.3\%. Figure \ref{fig:AmedeaConv} shows the convergence profiles obtained by GMRES(100) preconditioned with \precnameVar{1} for the saddle-point matrix $\mathcal{A}$ and AMG as implemented in \cite{IsoFriSpiJan21,paludetto2019novel} for the leading block $A$. The same AMG was used to approximately apply $\widetilde{S}_u^{-1}$ in the \precnameVar{1} algorithm. On the leading block $A$, the AMG algorithm has grid complexity 1.048
and operator complexity 1.323, while for $\widetilde{S}_u^{-1}$ they read 1.044
and 1.418, respectively. Hence, in practice there is no increase in the cost for using the same AMG algorithm to $A$ and $S_u$.
The iteration count for $A$ and $\mathcal{A}$ turns to be also the same, so that we can conclude that solving the problem with or without the faults implies practically the same computational cost. Notice also the very fast convergence rate (less than 30 iterations) for such a challenging test case.

\begin{figure}
  \centering
  \null\hfill
  \subfloat[]{\includegraphics[height=0.35\textwidth]{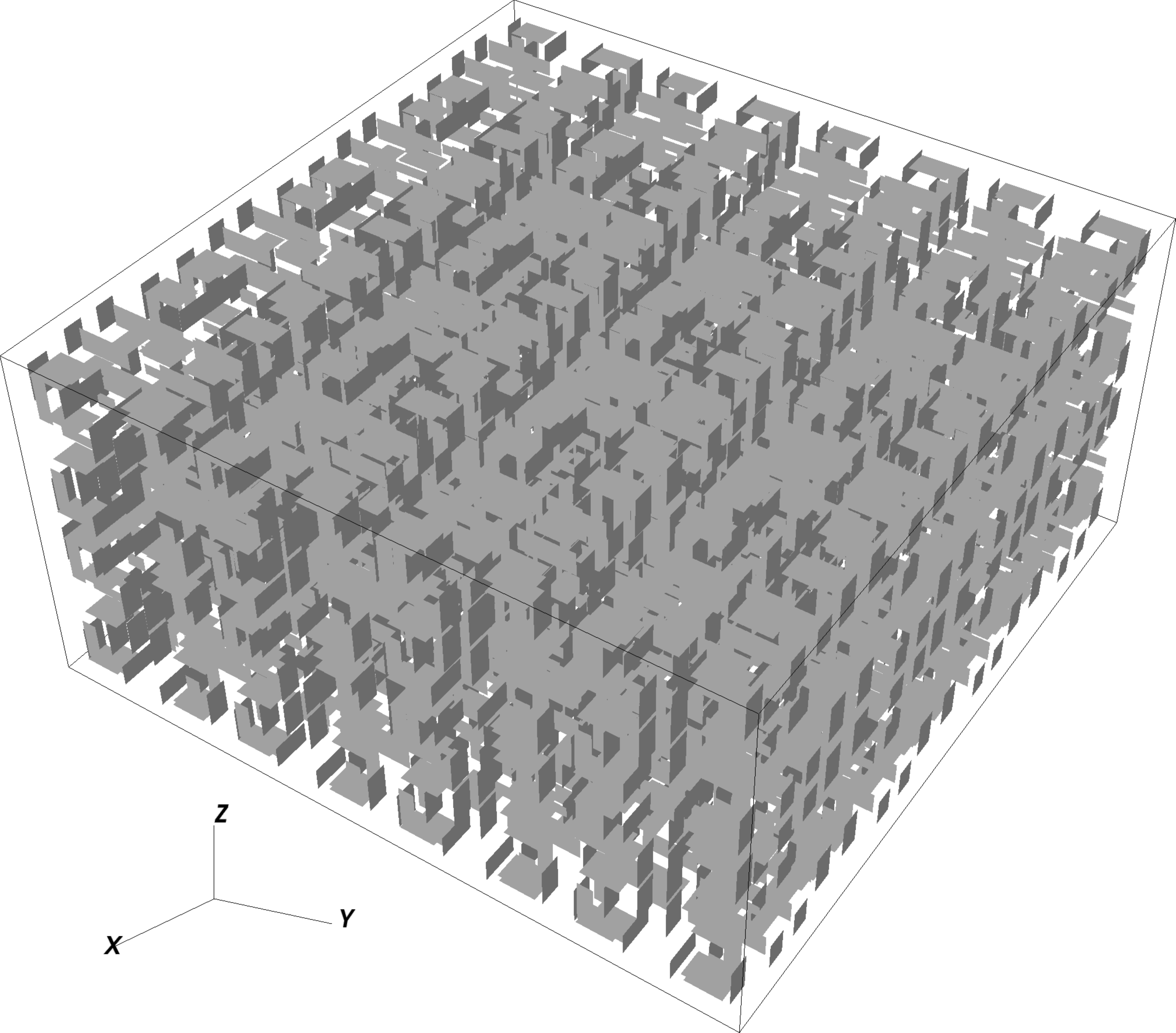}\label{fig:faultsMesh}}
  \hfill
  \subfloat[]{\includegraphics[width=0.50\textwidth]{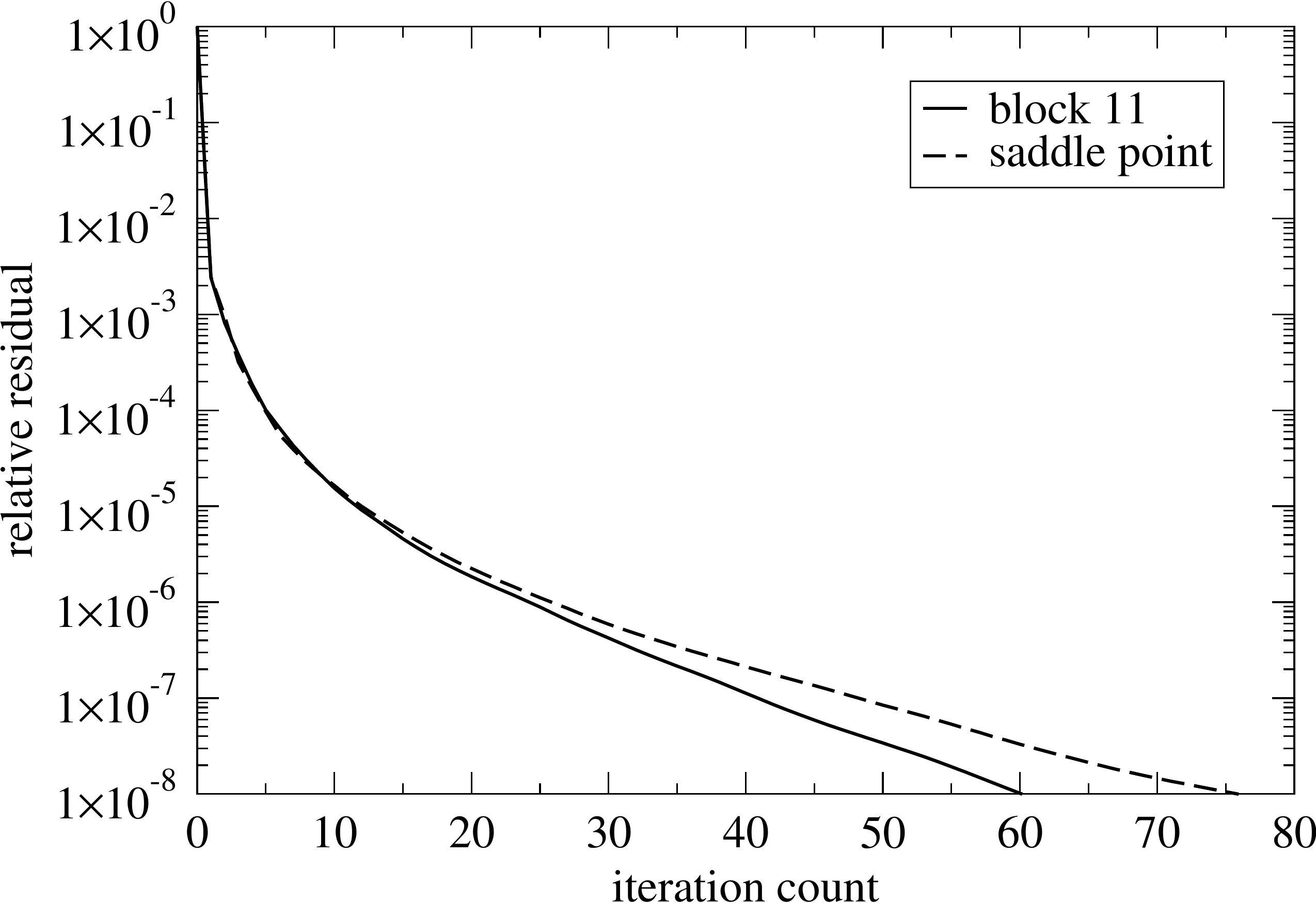}}
  \hfill\null
  \caption{Test B: Computational grid (a) and convergence profiles (b) for leading block and saddle
    point matrix.}
  \label{fig:faultsRes}
\end{figure}

As a last test case, denoted as {\tt Test B}, we consider a large-size problem also used to investigate the parallel efficiency of the proposed approach. 
The cube shown in Figure \ref{fig:faultsRes}, consisting of 95,451 nodes, 512,771 hexahedral elements and 29,324 interface elements, which correspond to 18 fractures, is repeatedly mirrored in order to generate a problem with an increasing size. For example,
mirroring the original cube
8 times gives rise to a
computational grid with 24,347,657 nodes, 129,736,704 elements and 7,506,944 interface
elements, i.e., 4,608 fractures. The problem sizes are $n_u = 73,042,971$ and $n_t = 10,553,856$, with
$\text{nnz}(A) = 3,218,602,977$ and $\text{nnz}(B) = 63,323,136$. The $n_t/n_u$ ratio is
14.4\%. As for {\tt Test A}, we report in Figure \ref{fig:faultsRes} the convergence profiles
for both the leading block $A$ and the saddle-point matrix $\mathcal{A}$.
The AMG algorithm used for $A$ and $S_u$ is the same, with grid and operator complexity equal to 1.053 and 1.055, and 1.363 and 1.393, respectively. 
In other words, the computational cost for applying the AMG preconditioner to $A$ and the overall \precnameVar{1} algorithm is roughly the same.
As shown in Figure \ref{fig:faultsRes}, the convergence is also very similar and quite fast, achieving the desired accuracy in less than 80 iterations for this huge problem.


Strong and weak scalability analyses are assessed to investigate the \precname suitability to
high performance computing. To this aim, we run a set of simulations on the Marconi100 supercomputer at the Italian Center for High Performance Computing (CINECA). Marconi100 is based on the IBM Power9 architecture with 980 nodes, each one equipped with two 16-core 
AC922 @3.1 GHz processors and 256 GB of RAM. 
The parallel implementation described in Section \ref{sec:impl} is used.
For the strong scalability test, we consider the {\tt Test B} configuration mentioned above and progressively increase the
number of processes, while in 
the weak scalability test we keep a constant work load
per process by properly varying the number of mirroring of the basic fractured cube. 
In all tests, we bind 4 openMP threads
to each MPI rank and vary the total number of MPI processes. This configuration, that is 4 threads per rank,
is the one typically giving the lowest solution time on Marconi100.

\begin{figure}
  \centering
  \null\hfill
  \includegraphics[height=0.28\textwidth]{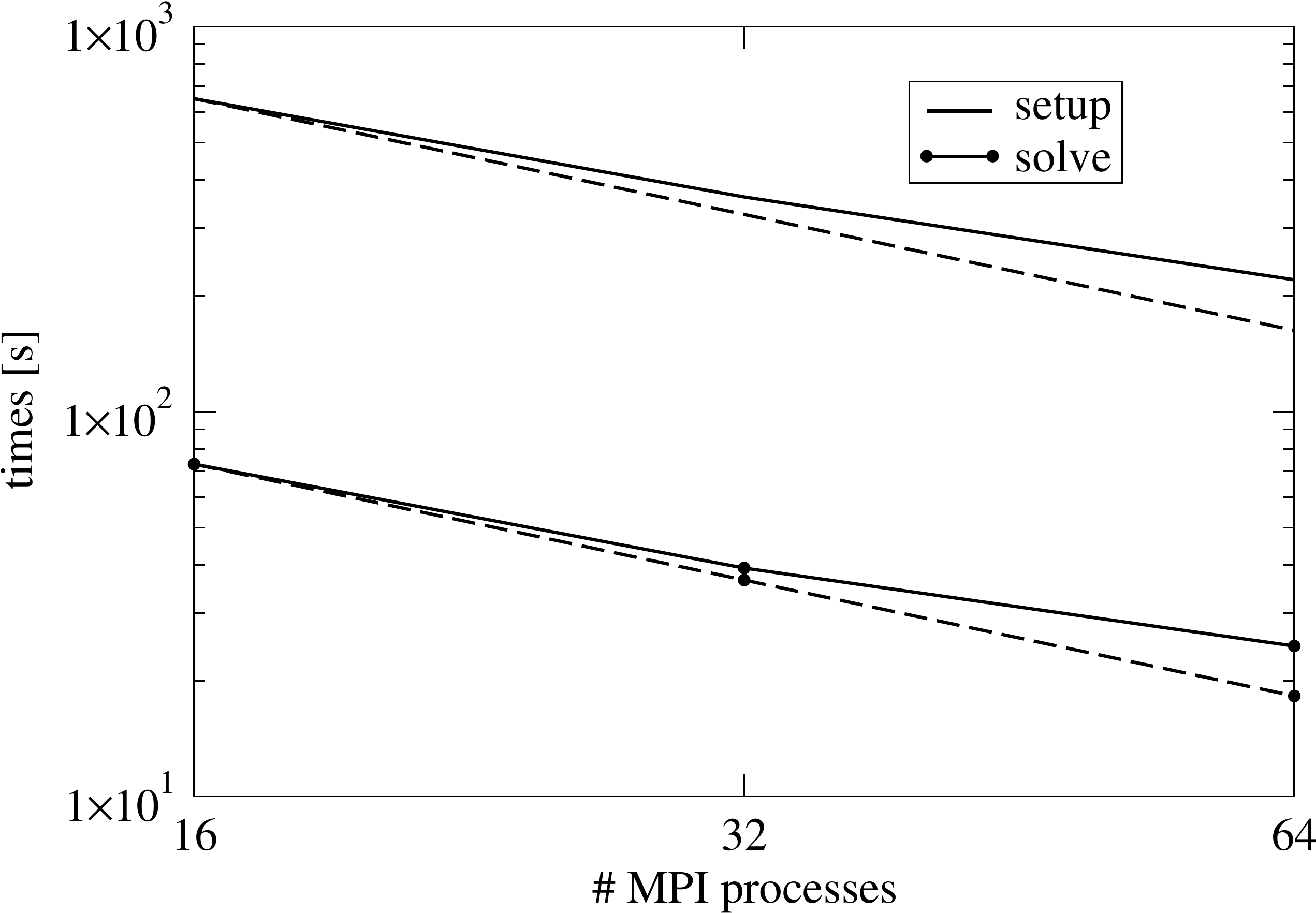}\hfill
  \includegraphics[height=0.30\textwidth]{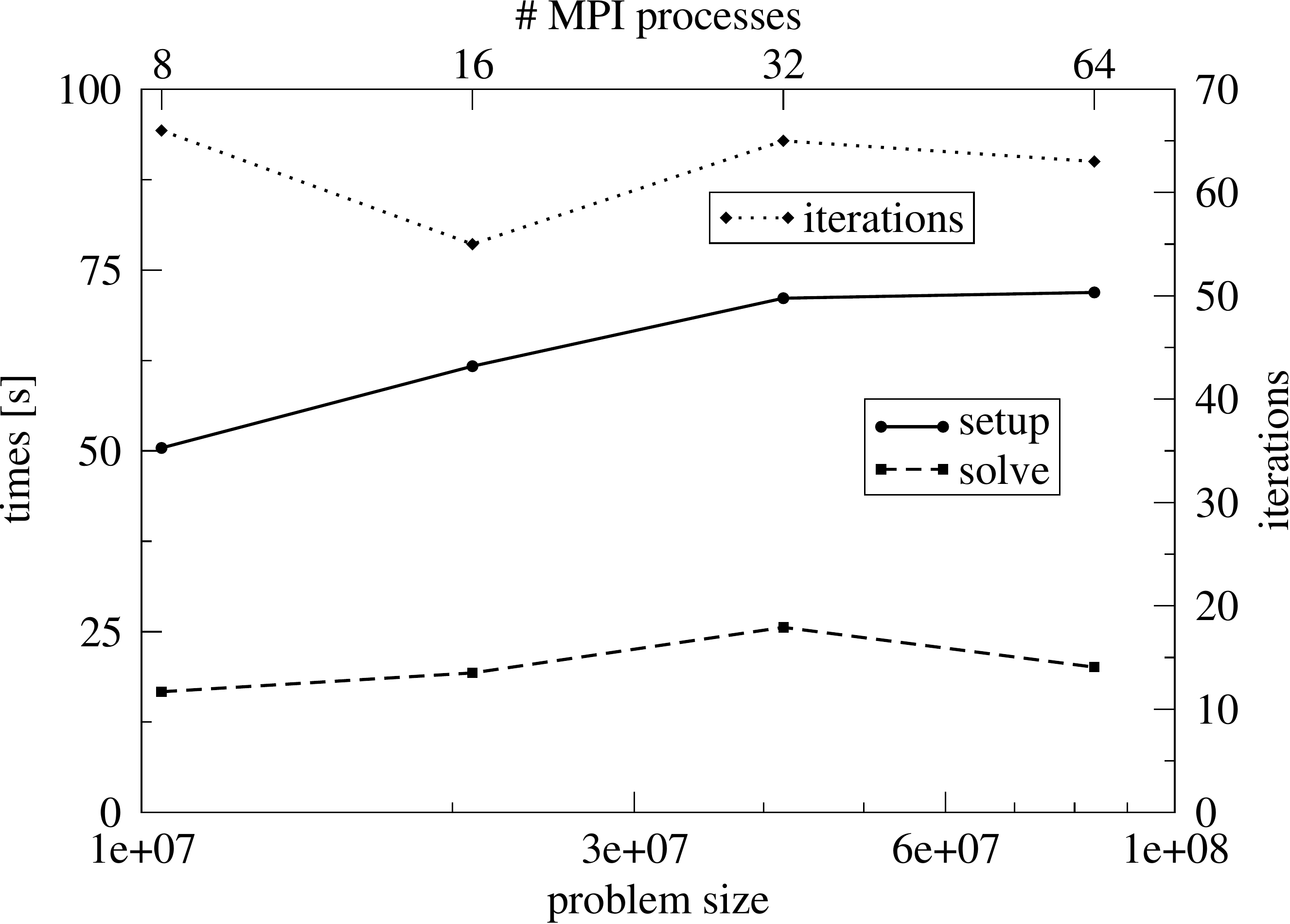}
  \hfill\null
  \caption{Test B: Strong (left) and weak (right) scalability tests.}
  \label{fig:faultsScal}
\end{figure}

Figure \ref{fig:faultsScal} shows the strong (left) and weak (right) scalability outcomes for both
setup and solving times. It can be noticed that both the strong and the weak scalability are quite close to the
ideal one. 
The only exception appears to be the setup time in the weak scalability, which
slightly increases from 50.4 s to 71.9 s while increasing the problem size by a factor
8. This suboptimal behavior can be ascribed to the extremely small computational load per process that
affects the lower AMG levels. In fact, independently of the actual problem size, 
the AMG hierarchy leads to very small coarse levels, that are difficult to handle efficiently
especially when significant computational resources are allocated. This inefficiency on the lower
levels is more pronounced for larger problems and it is why the weak scalability slightly strays from the ideal. 
Anyway, the overall parallel efficiency of the proposed algorithm, in both the strong and weak scalability tests, appears to be quite satisfactory.


\section{Conclusions}
\label{sec:concl}

A novel preconditioning approach is presented with the aim of accelerating
the solution of saddle-point systems arising from contact problems modeled with the aid of Lagrange multipliers. The proposed algorithm belongs to the class
of 
constraint preconditioners, but differs from the classical approaches with saddle-point systems because 
the primal Schur complement is approximated by introducing an appropriate augmentation to the leading block.
The properties of this strategy are investigated by addressing both mid- and large-size
problems, 
showing algorithmic robustness, computational efficiency and parallel scalability.
The main outcomes of the research can be summarized as follows.
\begin{itemize}
    \item Since the \precname algorithm belongs to the class of constraint preconditioners, a set of bounds for the eigenspectrum of the preconditioned matrix can be proved 
    as well. 
    Though these bounds often prove to be pretty loose, they can provide useful theoretical indications on the expected convergence behavior of the method.
    \item The reverse and augmented approach allows to address effectively problems with a singular leading block. 
    Even when the leading block is non-singular, the proposed augmentation strategy, which is based on the use of local pieces of information, generates a primal Schur complement that generally preserves, and in many instances also improves, the conditioning of the original leading block. This property is fundamental to keep under control the overall \precname computational cost.
    \item 
    As compared to the traditional constraint preconditioning used in saddle-point systems,
    the \precname set-up and application is generally cheaper, 
    more robust and more effective in accelerating the GMRES convergence.
    \item 
    \precname can be efficiently implemented and used in
    high performance computations, 
    provided that 
    effective kernels are available for 
    basic sparse matrix operations. In this work,
    the Chronos package has been used to provide the mid-level kernels, exhibiting a nearly optimal scalability up to
    hundreds of cores.
\end{itemize}

The algebraic properties of the \precname approach can be easily extended to other applications governed by saddle-point problems. 
The source files of the \precname implementation used in this work are available on GitHub at \url{https://github.com/matteofrigo5/HPC\_ReverseAugmentedConstrained}. 



\section*{Acknowledgements}
\label{sec::acknow}
Portions of this work were performed 
within the 2020 INdAM-GNCS project ``Optimization and advanced linear algebra for PDE-governed problems''.


\begin{thebibliography}{10}
\expandafter\ifx\csname url\endcsname\relax
  \def\url#1{\texttt{#1}}\fi
\expandafter\ifx\csname urlprefix\endcsname\relax\def\urlprefix{URL }\fi
\expandafter\ifx\csname href\endcsname\relax
  \def\href#1#2{#2} \def\path#1{#1}\fi

\bibitem{KikOde88}
N.~Kikuchi, J.~T. Oden, {Contact Problems in Elasticity: A Study of Variational
  Inequalities and Finite Element Methods}, SIAM, Philadelphia, PA, USA, 1988.
\newblock \href {https://doi.org/10.1137/1.9781611970845}
  {\path{doi:10.1137/1.9781611970845}}.

\bibitem{Lau03}
T.~A. Laursen, {Computational Contact and Impact Mechanics: Fundamentals of
  Modeling Interfacial Phenomena in Nonlinear Finite Element Analysis},
  Springer-Verlag Berlin Heidelberg, 2003.
\newblock \href {https://doi.org/10.1007/978-3-662-04864-1}
  {\path{doi:10.1007/978-3-662-04864-1}}.

\bibitem{Wri06}
P.~Wriggers, {Computational Contact Mechanics}, 2nd Edition, Springer-Verlag
  Berlin Heidelberg, 2006.
\newblock \href {https://doi.org/10.1007/978-3-540-32609-0}
  {\path{doi:10.1007/978-3-540-32609-0}}.

\bibitem{SimHug98}
J.~C. Simo, T.~J.~R. Hughes, {Computational Inelasticity}, Springer-Verlag New
  York, 1998.
\newblock \href {https://doi.org/10.1007/b98904} {\path{doi:10.1007/b98904}}.

\bibitem{Zie00}
O.~C. Zienkiewicz, R.~L. Taylor, {The Finite Element Method: Solid Mechanics},
  Vol.~2, Butterworth-Heinemann, 2000.

\bibitem{Bat06}
K.~J. Bathe, {Finite Element Procedures}, Klaus-Jurgen Bathe, 2006.

\bibitem{KhoLew99}
A.~R. Khoei, R.~W. Lewis, {Adaptive finite element remeshing in a large
  deformation analysis of metal powder forming}, Int. J. Numer. Meth. Eng.
  45~(7) (1999) 801--820.

\bibitem{OnaRoj04}
E.~Onate, J.~Rojek, {Combination of discrete element and finite element methods
  for dynamic analysis of geomechanics problems}, Comp. Meth. Appl. Mech. Eng.
  193~(27) (2004) 3087--3128.
\newblock \href {https://doi.org/10.1016/j.cma.2003.12.056}
  {\path{doi:10.1016/j.cma.2003.12.056}}.

\bibitem{FerGamJanTea08}
M.~Ferronato, G.~Gambolati, C.~Janna, P.~Teatini, {Numerical modelling of
  regional faults in land subsidence prediction above gas/oil reservoirs}, Int.
  J. Numer. Anal. Meth. Geomech. 32~(6) (2008) 633--657.
\newblock \href {https://doi.org/10.1002/nag.640} {\path{doi:10.1002/nag.640}}.

\bibitem{BenEssFak16}
E.~H. Benkhira, E.~H. Essoufi, R.~Fakhar, {On convergence of the penalty method
  for a static unilateral contact problem with nonlocal friction in
  electro-elasticity}, European J. Appl. Math. 27~(1) (2016) 1--22.
\newblock \href {https://doi.org/10.1017/S0956792515000248}
  {\path{doi:10.1017/S0956792515000248}}.

\bibitem{BurErn17}
E.~Burman, A.~Ern, {A nonlinear consistent penalty method for positivity
  preservation in the finite element approximation of the transport equation},
  Comp. Meth. Appl. Mech. Eng. 320 (2017) 122--132.
\newblock \href {https://doi.org/10.1016/j.cma.2017.03.019}
  {\path{doi:10.1016/j.cma.2017.03.019}}.

\bibitem{Ber84}
D.~P. Bertsekas, {Constrained Optimization and Lagrange Multiplier Methods},
  Academic Press New York, 1984.

\bibitem{HagHueWoh08}
C.~Hager, S.~H\"{u}eber, B.~I. Wohlmuth, {A stable energy-conserving approach
  for frictional contact problems based on quadrature formulas}, Int. J. Numer.
  Meth. Eng. 73~(2) (2008) 205--225.
\newblock \href {https://doi.org/10.1002/nme.2069}
  {\path{doi:10.1002/nme.2069}}.

\bibitem{fraferjantea16}
A.~Franceschini, M.~Ferronato, C.~Janna, P.~Teatini, {A novel Lagrangian
  approach for the stable numerical simulation of fault and fracture
  mechanics}, J. Comput. Phys. 314 (2016) 503--521.
\newblock \href {https://doi.org/10.1016/j.jcp.2016.03.032}
  {\path{doi:10.1016/j.jcp.2016.03.032}}.

\bibitem{benzi2005}
M.~Benzi, G.~H. Golub, J.~Liesen, {Numerical solution of saddle point
  problems}, Acta Numer. 14 (2005) 1--137.
\newblock \href {https://doi.org/10.1017/S0962492904000212}
  {\path{doi:10.1017/S0962492904000212}}.

\bibitem{li2001existence}
D.~Li, X.~L. Sun, {Existence of a saddle point in nonconvex constrained
  optimization}, J. Glob. Optim. 21~(1) (2001) 39--50.
\newblock \href {https://doi.org/10.1023/A:1017970111378}
  {\path{doi:10.1023/A:1017970111378}}.

\bibitem{Bergamaschi2007137}
L.~Bergamaschi, J.~Gondzio, M.~Venturin, G.~Zilli, {Inexact constraint
  preconditioners for linear systems arising in interior point methods},
  Comput. Optim. Appl. 36~(2-3) (2007) 137--147.
\newblock \href {https://doi.org/10.1007/s10589-006-9001-0}
  {\path{doi:10.1007/s10589-006-9001-0}}.

\bibitem{schoberl2007symmetric}
J.~Sch{\"o}berl, W.~Zulehner, {Symmetric indefinite preconditioners for saddle
  point problems with applications to PDE-constrained optimization problems},
  SIAM J. Matrix Anal. Appl. 29~(3) (2007) 752--773.
\newblock \href {https://doi.org/10.1137/060660977}
  {\path{doi:10.1137/060660977}}.

\bibitem{Pearson2020}
J.~W. Pearson, M.~Porcelli, M.~Stoll, {Interior-point methods and
  preconditioning for PDE-constrained optimization problems involving sparsity
  terms}, Numer. Lin. Alg. Appl. 27~(2) (2020).
\newblock \href {https://doi.org/10.1002/nla.2276}
  {\path{doi:10.1002/nla.2276}}.

\bibitem{vassilevski1996prec}
P.~S. Vassilevski, R.~D. Lazarov, {Preconditioning mixed finite element
  saddle-point elliptic problems}, Numer. Linear Algebra Appl. 3~(1) (1996)
  1--20.
\newblock \href
  {https://doi.org/10.1002/(SICI)1099-1506(199601/02)3:1<1::AID-NLA67>3.0.CO;2-E}
  {\path{doi:10.1002/(SICI)1099-1506(199601/02)3:1<1::AID-NLA67>3.0.CO;2-E}}.

\bibitem{barrientos2002mixed}
M.~A. Barrientos, G.~N. Gatica, E.~P. Stephan, {A mixed finite element method
  for nonlinear elasticity: two-fold saddle point approach and a-posteriori
  error estimate}, Numer. Math. 91~(2) (2002) 197--222.
\newblock \href {https://doi.org/10.1007/s002110100337}
  {\path{doi:10.1007/s002110100337}}.

\bibitem{loghin2004analysis}
D.~Loghin, A.~J. Wathen, {Analysis of preconditioners for saddle-point
  problems}, SIAM J. Sci. Comput. 25~(6) (2004) 2029--2049.
\newblock \href {https://doi.org/10.1137/S1064827502418203}
  {\path{doi:10.1137/S1064827502418203}}.

\bibitem{olshanskii2006uniform}
M.~A. Olshanskii, J.~Peters, A.~Reusken, {Uniform preconditioners for a
  parameter dependent saddle point problem with application to generalized
  Stokes interface equations}, Numer. Math. 105~(1) (2006) 159--191.
\newblock \href {https://doi.org/10.1007/s00211-006-0031-4}
  {\path{doi:10.1007/s00211-006-0031-4}}.

\bibitem{bergamaschi2008mixed}
L.~Bergamaschi, M.~Ferronato, G.~Gambolati, {Mixed constraint preconditioners
  for the iterative solution of FE coupled consolidation equations}, J. Comput.
  Phys. 227~(23) (2008) 9885--9897.
\newblock \href {https://doi.org/10.1016/j.jcp.2008.08.002}
  {\path{doi:10.1016/j.jcp.2008.08.002}}.

\bibitem{axelsson2012stable}
O.~Axelsson, R.~Blaheta, P.~Byczanski, {Stable discretization of poroelasticity
  problems and efficient preconditioners for arising saddle point type
  matrices}, Comput. Visual. Sci. 15~(4) (2012) 191--207.
\newblock \href {https://doi.org/10.1007/s00791-013-0209-0}
  {\path{doi:10.1007/s00791-013-0209-0}}.

\bibitem{cao2015relaxed}
Y.~Cao, J.-L. Dong, Y.-M. Wang, A relaxed deteriorated pss preconditioner for
  nonsymmetric saddle point problems from the steady navier--stokes equation,
  J. Comput. Appl. Math. 273 (2015) 41--60.
\newblock \href {https://doi.org/10.1016/j.cam.2014.06.001}
  {\path{doi:10.1016/j.cam.2014.06.001}}.

\bibitem{Pearson2018331}
J.~W. Pearson, J.~Pestana, D.~J. Silvester, {Refined saddle-point
  preconditioners for discretized Stokes problems}, Numer. Math. 138~(2) (2018)
  331--363.
\newblock \href {https://doi.org/10.1007/s00211-017-0908-4}
  {\path{doi:10.1007/s00211-017-0908-4}}.

\bibitem{hong2020parameter}
Q.~Hong, J.~Kraus, M.~Lymbery, F.~Philo, {Parameter-robust Uzawa-type iterative
  methods for double saddle point problems arising in Biot's consolidation and
  multiple-network poroelasticity models}, Math. Models Methods Appl. Sci.
  30~(13) (2020) 2523--2555.
\newblock \href {https://doi.org/10.1142/S0218202520500499}
  {\path{doi:10.1142/S0218202520500499}}.

\bibitem{Murphy20001969}
M.~F. Murphy, G.~H. Golub, A.~J. Wathen, {Note on preconditioning for
  indefinite linear systems}, SIAM J. Sci. Comp. 21~(6) (2000) 1969--1972.
\newblock \href {https://doi.org/10.1137/S1064827599355153}
  {\path{doi:10.1137/S1064827599355153}}.

\bibitem{Keller20001300}
C.~Keller, N.~I.~M. Gould, A.~J. Wathen, {Constraint preconditioning for
  indefinite linear systems}, SIAM J. Matrix Anal. Appl. 21~(4) (2000)
  1300--1317.
\newblock \href {https://doi.org/10.1137/S0895479899351805}
  {\path{doi:10.1137/S0895479899351805}}.

\bibitem{Bergamaschi2004149}
L.~Bergamaschi, J.~Gondzio, G.~Zilli, {Preconditioning indefinite systems in
  interior point methods for optimization}, Comput. Optim. Appl. 28~(2) (2004)
  149--171.
\newblock \href {https://doi.org/10.1023/B:COAP.0000026882.34332.1b}
  {\path{doi:10.1023/B:COAP.0000026882.34332.1b}}.

\bibitem{Bergamaschi20072647}
L.~Bergamaschi, M.~Ferronato, G.~Gambolati, {Novel preconditioners for the
  iterative solution to FE-discretized coupled consolidation equations}, Comp.
  Meth. Appl. Mech. Eng. 196~(25-28) (2007) 2647--2656.
\newblock \href {https://doi.org/10.1016/j.cma.2007.01.013}
  {\path{doi:10.1016/j.cma.2007.01.013}}.

\bibitem{Janna2012661}
C.~Janna, M.~Ferronato, G.~Gambolati, {Parallel inexact constraint
  preconditioning for ill-conditioned consolidation problems}, Comput. Geosci.
  16~(3) (2012) 661--675.
\newblock \href {https://doi.org/10.1007/s10596-012-9276-4}
  {\path{doi:10.1007/s10596-012-9276-4}}.

\bibitem{ferronato2019general}
M.~Ferronato, A.~Franceschini, C.~Janna, N.~Castelletto, H.~A. Tchelepi, {A
  general preconditioning framework for coupled multiphysics problems with
  application to contact-and poro-mechanics}, J. Comput. Phys. 398 (2019)
  108887.
\newblock \href {https://doi.org/10.1016/j.jcp.2019.108887}
  {\path{doi:10.1016/j.jcp.2019.108887}}.

\bibitem{Nardean2021}
S.~Nardean, M.~Ferronato, A.~S. Abushaikha, {A novel block non-symmetric
  preconditioner for mixed-hybrid finite-element-based Darcy flow simulations},
  J. Comput. Phys. 442 (2021).
\newblock \href {https://doi.org/10.1016/j.jcp.2021.110513}
  {\path{doi:10.1016/j.jcp.2021.110513}}.

\bibitem{Silvester2001261}
D.~Silvester, H.~Elman, D.~Kay, A.~Wathen, {Efficient proconditioning of the
  linearized Navier-Stokes equations for incompressible flow}, J. Comput. Appl.
  Math. 128~(1-2) (2001) 261--279.
\newblock \href {https://doi.org/10.1016/S0377-0427(00)00515-X}
  {\path{doi:10.1016/S0377-0427(00)00515-X}}.

\bibitem{elman2002preconditioners}
H.~C. Elman, {Preconditioners for saddle point problems arising in
  computational fluid dynamics}, Appl. Numer. Math. 43~(1-2) (2002) 75--89.
\newblock \href {https://doi.org/10.1016/S0168-9274(02)00118-6}
  {\path{doi:10.1016/S0168-9274(02)00118-6}}.

\bibitem{cao2003fast}
Z.-H. Cao, Fast uzawa algorithm for generalized saddle point problems, Appl.
  Numer. Math. 46~(2) (2003) 157--171.
\newblock \href {https://doi.org/10.1016/S0168-9274(03)00023-0}
  {\path{doi:10.1016/S0168-9274(03)00023-0}}.

\bibitem{Elman20061651}
H.~Elman, V.~E. Howle, J.~Shadid, R.~Shuttleworth, R.~Tuminaro, {Block
  preconditioners based on approximate commutators}, SIAM J. Sci. Comp. 27~(5)
  (2006) 1651--1668.
\newblock \href {https://doi.org/10.1137/040608817}
  {\path{doi:10.1137/040608817}}.

\bibitem{choi2015practical}
Y.~Choi, C.~Farhat, W.~Murray, M.~Saunders, {A practical factorization of a
  Schur complement for PDE-constrained distributed optimal control}, J. Sci.
  Comput. 65~(2) (2015) 576--597.
\newblock \href {https://doi.org/10.1007/s10915-014-9976-0}
  {\path{doi:10.1007/s10915-014-9976-0}}.

\bibitem{Castelletto2016894}
N.~Castelletto, J.~A. White, M.~Ferronato, {Scalable algorithms for three-field
  mixed finite element coupled poromechanics}, J. Comput. Phys. 327 (2016)
  894--918.
\newblock \href {https://doi.org/10.1016/j.jcp.2016.09.063}
  {\path{doi:10.1016/j.jcp.2016.09.063}}.

\bibitem{Greif20191}
C.~Greif, M.~Wathen, {Conjugate gradient for nonsingular saddle-point systems
  with a maximally rank-deficient leading block}, J. Comput. Appl. Math. 358
  (2019) 1--11.
\newblock \href {https://doi.org/10.1016/j.cam.2019.02.016}
  {\path{doi:10.1016/j.cam.2019.02.016}}.

\bibitem{fortin1983augmented}
M.~Fortin, R.~Glowinski, {Augmented Lagrangian methods: applications to the
  numerical solution of boundary-value problems}, Elsevier, 1983.
\newblock \href {https://doi.org/10.1016/s0168-2024(08)x7003-1}
  {\path{doi:10.1016/s0168-2024(08)x7003-1}}.

\bibitem{bacuta2006unified}
C.~Bacuta, {A unified approach for Uzawa algorithms}, SIAM J. Numer. Anal.
  44~(6) (2006) 2633--2649.
\newblock \href {https://doi.org/10.1137/050630714}
  {\path{doi:10.1137/050630714}}.

\bibitem{Benzi20062095}
M.~Benzi, M.~A. Olshanskii, {An augmented Lagrangian-based approach to the
  Oseen problem}, SIAM J. Sci. Comp. 28~(6) (2006) 2095--2113.
\newblock \href {https://doi.org/10.1137/050646421}
  {\path{doi:10.1137/050646421}}.

\bibitem{lee2007robust}
Y.-J. Lee, J.~Wu, J.~Xu, L.~Zikatanov, {Robust subspace correction methods for
  nearly singular systems}, Math. Models Methods Appl. Sci. 17~(11) (2007)
  1937--1963.
\newblock \href {https://doi.org/10.1142/S0218202507002522}
  {\path{doi:10.1142/S0218202507002522}}.

\bibitem{brezzi2012mixed}
F.~Brezzi, M.~Fortin, {Mixed and hybrid finite element methods}, Vol.~15,
  Springer Science \& Business Media, 2012.

\bibitem{CHRONOS-webpage}
M.~Frigo, G.~Isotton, C.~Janna, \href{https://www.m3eweb.it/chronos}{{Chronos}
  {W}eb page}, \url{https://www.m3eweb.it/chronos} (2020).
\newline\urlprefix\url{https://www.m3eweb.it/chronos}

\bibitem{franceschini2019block}
A.~Franceschini, N.~Castelletto, M.~Ferronato, {Block preconditioning for
  fault/fracture mechanics saddle-point problems}, Comput. Meth. Appl. Mech.
  Eng. 344 (2019) 376--401.
\newblock \href {https://doi.org/10.1016/j.cma.2018.09.039}
  {\path{doi:10.1016/j.cma.2018.09.039}}.

\bibitem{bergamaschi2012eigenvalue}
L.~Bergamaschi, {On eigenvalue distribution of constraint-preconditioned
  symmetric saddle point matrices}, Numer. Linear Algebra Appl. 19~(4) (2012)
  754--772.
\newblock \href {https://doi.org/10.1002/nla.806} {\path{doi:10.1002/nla.806}}.

\bibitem{axelsson2006eigenvalue}
O.~Axelsson, M.~Neytcheva, {Eigenvalue estimates for preconditioned saddle
  point matrices}, Numer. Linear Algebra Appl. 13~(4) (2006) 339--360.
\newblock \href {https://doi.org/10.1002/nla.469} {\path{doi:10.1002/nla.469}}.

\bibitem{ruiz2018refining}
D.~Ruiz, A.~Sartenaer, C.~Tannier, {Refining the lower bound on the positive
  eigenvalues of saddle point matrices with insights on the interactions
  between the blocks}, SIAM J. Matrix Anal. Appl. 39~(2) (2018) 712--736.
\newblock \href {https://doi.org/10.1137/16M108152X}
  {\path{doi:10.1137/16M108152X}}.

\bibitem{Aagaard20133059}
B.~Aagaard, M.~Knepley, C.~Williams, A domain decomposition approach to
  implementing fault slip in finite-element models of quasi-static and dynamic
  crustal deformation, J. Geophys. Res.: Solid Earth 118~(6) (2013) 3059--3079.
\newblock \href {https://doi.org/10.1002/jgrb.50217}
  {\path{doi:10.1002/jgrb.50217}}.

\bibitem{Jha20143776}
B.~Jha, R.~Juanes, Coupled multiphase flow and poromechanics: A computational
  model of pore pressure effects on fault slip and earthquake triggering, Water
  Resour. Res. 50~(5) (2014) 3776--3808.
\newblock \href {https://doi.org/10.1002/2013WR015175}
  {\path{doi:10.1002/2013WR015175}}.

\bibitem{White201655}
J.~White, N.~Castelletto, H.~Tchelepi, Block-partitioned solvers for coupled
  poromechanics: A unified framework, Comp. Meth. Appl. Mech. Eng. 303 (2016)
  55--74.
\newblock \href {https://doi.org/10.1016/j.cma.2016.01.008}
  {\path{doi:10.1016/j.cma.2016.01.008}}.

\bibitem{Pow69}
M.~J.~D. Powell, {A Method for Nonlinear Constraints in Minimization Problems},
  Academic Press London, 1969.

\bibitem{Golub20032076}
G.~Golub, C.~Greif, On solving block-structured indefinite linear systems, SIAM
  J. Sci. Comp. 24~(6) (2003) 2076--2092.
\newblock \href {https://doi.org/10.1137/S1064827500375096}
  {\path{doi:10.1137/S1064827500375096}}.

\bibitem{saad1986gmres}
Y.~Saad, M.~H. Schultz, {GMRES: A generalized minimal residual algorithm for
  solving nonsymmetric linear systems}, SIAM J. Sci. Stat. Comput. 7~(3) (1986)
  856--869.
\newblock \href {https://doi.org/10.1137/0907058} {\path{doi:10.1137/0907058}}.

\bibitem{IsoFriSpiJan21}
G.~Isotton, M.~Frigo, N.~Spiezia, C.~Janna, {Chronos: A general purpose
  classical AMG solver for High Performance Computing}, SIAM J. Sci. Comput.
  43~(5) (2021) C335--C357.
\newblock \href {https://doi.org/10.1137/21M1398586}
  {\path{doi:10.1137/21M1398586}}.

\bibitem{pMETIS}
G.~Karypis, K.~Schloegel,
  \href{http://glaros.dtc.umn.edu/gkhome/fetch/sw/parmetis/manual.pdf}{{Parallel
  Graph Partitioning and Sparse Matrix Ordering Library, Version 4.0}} (2013).
\newline\urlprefix\url{http://glaros.dtc.umn.edu/gkhome/fetch/sw/parmetis/manual.pdf}

\bibitem{franceschini2015modelling}
A.~Franceschini, P.~Teatini, C.~Janna, M.~Ferronato, G.~Gambolati, S.~Ye,
  D.~Carre{\'o}n-Freyre, {Modelling ground rupture due to groundwater
  withdrawal: applications to test cases in China and Mexico}, Proc. Int.
  Assoc. Hydrol. 372 (2015) 63--68.
\newblock \href {https://doi.org/10.5194/piahs-372-63-2015}
  {\path{doi:10.5194/piahs-372-63-2015}}.

\bibitem{paludetto2019novel}
V.~A. Paludetto~Magri, A.~Franceschini, C.~Janna, {A novel algebraic multigrid
  approach based on adaptive smoothing and prolongation for ill-conditioned
  systems}, SIAM J. Sci. Comput. 41~(1) (2019) A190--A219.
\newblock \href {https://doi.org/10.1137/17M1161178}
  {\path{doi:10.1137/17M1161178}}.

\bibitem{janna2015fsaipack}
C.~Janna, M.~Ferronato, F.~Sartoretto, G.~Gambolati, {FSAIPACK: A software
  package for high-performance factored sparse approximate inverse
  preconditioning}, ACM Trans. Math. Software (TOMS) 41~(2) (2015) 1--26.
\newblock \href {https://doi.org/10.1145/2629475} {\path{doi:10.1145/2629475}}.

\bibitem{franceschini2019robust}
A.~Franceschini, V.~A.~P. Magri, G.~Mazzucco, N.~Spiezia, C.~Janna, {A robust
  adaptive algebraic multigrid linear solver for structural mechanics}, Comput.
  Meth. Appl. Mech. Eng. 352 (2019) 389--416.
\newblock \href {https://doi.org/10.1016/j.cma.2019.04.034}
  {\path{doi:10.1016/j.cma.2019.04.034}}.

\bibitem{castelletto2013geological}
N.~Castelletto, G.~Gambolati, P.~Teatini, {Geological CO$_2$ sequestration in
  multi-compartment reservoirs: Geomechanical challenges}, J. Geophys. Res.
  Solid Earth 118~(5) (2013) 2417--2428.
\newblock \href {https://doi.org/10.1002/jgrb.50180}
  {\path{doi:10.1002/jgrb.50180}}.

\bibitem{castelletto2013multiphysics}
N.~Castelletto, P.~Teatini, G.~Gambolati, D.~Bossie-Codreanu, O.~Vinck{\'e},
  J.-M. Daniel, A.~Battistelli, M.~Marcolini, F.~Donda, V.~Volpi, {Multiphysics
  modeling of CO$_2$ sequestration in a faulted saline formation in Italy},
  Adv. Water Resour. 62 (2013) 570--587.
\newblock \href {https://doi.org/10.1016/j.advwatres.2013.04.006}
  {\path{doi:10.1016/j.advwatres.2013.04.006}}.

\end{thebibliography}

\end{document}